\documentclass[12pt]{article}
\usepackage{amsmath}
\usepackage{amssymb}
\usepackage{amsfonts}
\usepackage{eucal}
\usepackage[usenames]{color}
\usepackage{graphicx}
\numberwithin{equation}{section}
\newtheorem{Theorem}{Theorem}[section]
\newtheorem{Proposition}[Theorem]{Proposition}

\newtheorem{Remark}[Theorem]{Remark}

\begin{document}
\title{Optimization in the  first-passage problem of a diffusion  with Poissonian resetting}
 \author{Mario Abundo\thanks{Dipartimento di Matematica, Universit\`a
``Tor Vergata'', via della Ricerca Scientifica, I-00133 Roma,
Italy. E-mail: \texttt{abundo@mat.uniroma2.it}}  }

\date{}
\maketitle

\begin{abstract}
\noindent
We investigate the first-passage time (FPT) and first-exit time (FET) of a one-dimensi--
\par\noindent
onal, time-homogeneous diffusion process subject to Poissonian resetting. We first derive a general analytical relationship that expresses the Laplace transform (LT) and the expected value of the FPT (and FET) for the process with resetting in terms of the LT of the FPT (and FET) of the underlying diffusion without resetting. This framework is then applied to determine the optimal resetting rate $r$ that minimizes the expected FPT (and FET). We provide explicit results for drifted Brownian motion and  Ornstein-Uhlenbeck (OU) process. For Brownian motion, we extend existing literature by considering the case where the initial position $x$ differs from the resetting position $x_R$, providing a comprehensive parametric analysis. For the OU process, we provide new insights into the minimization of the expected FPT. Our results demonstrate how a strategic choice of the resetting rate can effectively regularize and accelerate the passage through one or two boundaries.

\end{abstract}


\section{Introduction}
In this paper, we study the first-passage time (FPT) and first-exit time (FET) of a one-dimensional diffusion process with Poissonian resetting, denoted by ${\cal X}(t)$. This process is derived from a time-homogeneous diffusion $X(t)$. Our objectives are twofold: first, we establish a relationship between the Laplace transform (LT) of the FPT (and FET) of ${\cal X}(t)$ and the LT of the  FPT (and  FET) of the underlying process $X(t)$. From this relationship, we derive the expected FPT (and FET) of ${\cal X}(t)$ in terms of the LT of the FPT (and FET) of the original process.
In fact, directly deriving the LT and the mean of FPT (and FET) for ${\cal X}(t)$ typically requires solving  differential systems which are significantly more complex than in the no-resetting case
(see \cite{abundo:MCAP2025}, \cite{abundo:FPA2023}). \par
Second, we address the optimization problem of minimizing the expected FPT (and FET) of ${\cal X}(t)$ with respect to the resetting rate $r$.  \par
For processes like Brownian motion (BM), the expected FPT through a fixed threshold is known to be infinite. However, the introduction of resetting  renders the expected FPT  finite
(see e.g.\cite{abundo:FPA2023}). In the two-boundary case, specifically the expected FET of BM from an interval $(0,b),$ it is finite for any $r \ge 0$ (see \cite{abundo:MCAP2025}); notably, a judicious choice of the resetting rate $r$ can expedite the passage through the boundaries, effectively reducing the expected FET.
Because resetting can both regularize (make finite) and accelerate passage times, it is of significant interest to identify the optimal resetting rate $r$ that minimizes the expected FPT (or FET) for general diffusion processes.
\par
The process ${\cal X}(t)$ is constructed as follows.\par\noindent
Let $X(t)$ be a one-dimensional time-homogeneous diffusion process satisfying the stochastic differential equation (SDE):
\begin{equation} \label{diffueq}
dX(t)= \mu(X(t)) dt + \sigma(X(t)) dB_t,
\end{equation}
with initial condition $X(0)=x,$
where $B_t$ is a standard BM. We assume the drift $\mu(\cdot)$ and diffusion coefficient $\sigma(\cdot)$ satisfy the standard conditions for the existence and uniqueness of a strong solution. Resetting events occur according to a homogeneous Poisson process with rate $r > 0$. Between resets, ${\cal X}(t)$ evolves according to \eqref{diffueq}; upon a reset, the process instantly returns to a fixed position $x_R$ and starts afresh. Consequently, the inter-resetting times are independent and exponentially distributed with parameter $r$.
In other words, in any time interval $(t, t+ \Delta t),$ with $\Delta t \rightarrow 0 ^+, $ the process can pass from ${\cal X}(t)$  to the position $x_R$ with probability $r \Delta t  + o( \Delta t),$ or it can continue its evolution as $X(t)$ with probability $1- r \Delta t + o( \Delta t ).$ The process ${\cal X}(t)$ so obtained is called diffusion process with Poissonian resetting, or simply diffusion with resetting;
this model shares similarities with the process in \cite{dicre:03}, which studied M/M/1 queues with catastrophes and their continuous approximation as a Wiener process with jumps to state $0$.
\par
The study of the first-passage time of a diffusion processes, with or without resetting, has diverse applications. These include biological models for neuronal activity (see e.g. \cite{lanska:89}, \cite{norisa:85}), credit risk modeling in mathematical finance (see e.g. \cite{jackson:stapro09}), and queuing theory (see e.g. the discussion in \cite{abundo:stapro12}). While the FPT (and FET) for BM with resetting have been explored in recent literature (\cite{abundo:MCAP2025}, \cite{abundo:TPMS2024}, \cite{abundo:FPA2023}), this paper generalizes the inquiry to a broader class of diffusion processes. We provide explicit examples for drifted BM and  Ornstein-Uhlenbeck (OU) process, with resetting. For drifted BM, we analyze both one-and two-boundary cases. For the OU process, we focus on the single-boundary FPT due to the analytical complexity of the two-boundary LT.
While the minimization of the expected FPT (and FET) for BM with resetting has been studied for the case $x = x_R$ in \cite{evans:11}, \cite{Pal:19}, we extend this analysis to $x \neq x_R$. We provide a detailed parametric analysis with graphical and numerical results to identify the optimal $r$. To our knowledge, the OU process with resetting has not been previously studied in this context.
\par
The paper is organized as follows: Section 2 presents general results and the derivation of expected FPT/FET via LTs. Sections 3 and 4 provide explicit calculations of the expected FPT and FET, respectively, for specific diffusion models. Section 5 is dedicated to the minimization of these expected times with respect to $r$, and Section 6 offers concluding remarks.

\section{General results}
For any $C^2$ function $u(x),$ the infinitesimal generator of the diffusion with resetting ${\cal X}(t)$ is given by (see e.g. \cite{abundo:FPA2023}):
\begin{equation} \label{generator}
{\cal L}u(x) = \frac 1 2 \sigma ^2(x) u''(x) + \mu (x) u'(x) +r (u(x_R) -u(x)) \equiv L u(x) +r (u(x_R) -u(x)) ,
\end{equation}
where $u'(x)$ and $u''(x)$ stand for the first and second derivative of $u(x).$
Here, $Lu(x)= \frac 1 2 \sigma ^2(x) u''(x) + \mu (x) u'(x)$ represents the  ``diffusion part'' of the generator, i.e. that concerning the diffusion process $X(t).$ \par
For an initial position ${\cal X}(0)= x>0,$ and non-negative reset position $x_R ,$  let
\begin{equation}
\tau _0(x,r) = \inf \{t>0: {\cal X} (t) = 0 \ | \ {\cal X} (0)=X(0) =x  \}
\end{equation}
be the first-passage time (FPT) of ${\cal X} (t)$ through zero, under the condition that ${\cal X} (0)=x $
(the notation includes the dependence on $x$ and $r,$ but not on $x_R,$ for the sake of simplicity).
\par
Moreover, for $b >0$ and $x_R, x  \in (0,b),$ we denote by  $\tau _{0b} (x,r)$
the first-exit time (FET) of $\mathcal X  (t)$ from the interval $(0,b),$ under the condition that $\mathcal X(0)=x ,$ namely:
\begin{equation}
\tau _{0b}  (x,r)= \min \{t>0: \mathcal X(t) \notin (0,b) | \mathcal X(0) =x  \}.
\end{equation}
(once again, the notation does not include the dependence on $x_R.$) \par\noindent
Note that we consider an interval $(0,b)$  for simplicity, but our study can be easily extended to any interval $(a,b).$ \par
In the following, we will use the same notation,  $\tau (x,r),$ for both the FPT and the FET;
from the context it will be clear which case we are referring to. \par\noindent

\subsection{The Laplace transforms of the survival function and of the first-passage time}
Let $Q _r (x,t;x_R)= P[ \tau (x,r)>t]$ be the survival function of the FPT (or FET) $\tau (x,r),$ where we use the notation of \cite{guoyan:24}, \cite{huang:24}, which makes explicit the dependence on $r$ and $x_R.$
Then, $Q  _r(x,t;x_R)$ satisfies the backward equation (see e.g. \cite{guoyan:24}, \cite{huang:24}):
$$
\frac {\partial} {\partial t} Q  _r(x,t; x_R) = {\cal L} Q   _r(x,t;x_R)
$$
\begin{equation} \label{eqsurvivalprob}
= \frac 1 2 \sigma ^2 (x) \frac {\partial ^2} {\partial x^2 }Q  _r(x,t; x_R) + \mu (x) \frac {\partial} {\partial x }Q  _r(x,t; x_R) +r \left [Q  _r(x_R,t; x_R) -Q _r(x,t; x_R) \right ]
\end{equation}
with the initial condition
 \begin{equation}
 Q _r(x,0;x_R) = 1 ,
 \end{equation}
and the boundary conditions:
\begin{equation}
 Q  _r(0,t; x_R) =0,
 \end{equation}
in the case of the FPT (one boundary),
\begin{equation}
 Q  _r(0,t; x_R) = Q _r(b,t; x_R)=0,
 \end{equation}
 in the case of the FET (two boundaries). \par\noindent
The operator ${\cal L},$ defined by  \eqref{generator}, i.e. the infinitesimal generator of the diffusion  with resetting $\mathcal X(t),$ is meant as an operator which acts on  $Q  _r(x,t; x_R),$ as a function of $x.$\par\noindent
For $\lambda >0,$ let $\widehat Q _r(x, \lambda; x_R)= \int_0^ \infty e^ {-\lambda t } Q  _r(x,t; x_R) dt$ be the Laplace transform (LT) of
$Q _r(x,t;x_R).$
By taking the Laplace transform in both sides of \eqref{eqsurvivalprob}, we get that $\widehat Q _r(x, \lambda; x_R)$
satisfies the equation:
 \begin{equation} \label{eqLTsurvivalprob}
\frac 1 2 \sigma ^2(x) \frac {\partial ^2} {\partial x^2 }\widehat Q_r(x, \lambda; x_R) + \mu (x)  \frac {\partial} {\partial x } \widehat Q _r(x, \lambda; x_R) - (r + \lambda) \widehat Q _r(x, \lambda; x_R) + r \widehat Q _r(x_R, \lambda; x_R) = -1,
 \end{equation}
 subject to the boundary conditions
  \begin{equation}
  \widehat Q _r(0, \lambda; x_R) =0,
  \end{equation}
  for one boundary, and
    \begin{equation}
\widehat Q _r(0, \lambda; x_R) = \widehat Q _r(b, \lambda; x_R)=0,
 \end{equation}
 for two boundaries. \par
By using a renewal argument (see e.g. \cite{guoyan:24}), one obtains:
\begin{equation} \label{renewalforQ}
Q _r(x,t;x_R) = e^ {-rt}Q  _0(x,t)+r \int _0 ^t ds e^ {-rs }Q  _0(x_R,s)Q _r(x,t-s;x_R),
\end{equation}
where $Q_0(x,t)$ is the survival function of the FPT (or FET)  for the process without resetting $X(t),$ that is, for $r=0.$ \par\noindent
For $\lambda >0,$ let
$\widehat Q _0(x, \lambda )= \int_0^ \infty e^ {-\lambda t } Q _0(x,t) dt$ be the LT of
$Q _0(x,t).$
Then from \eqref{renewalforQ} one gets:
\begin{equation} \label{renewalforQ2}
\widehat Q _r(x, \lambda; x_R) = \frac {\widehat Q _0(x,\lambda +r)} {1-r\widehat Q _0(x_R,\lambda +r) }.
\end{equation}
\par
We denote by $M _r (x, \lambda; x_R)$ the LT of $\tau (x,r),$ that is,
\begin{equation}
 M _r (x, \lambda; x_R) =  E \left [e^ {- \lambda  \tau (x,r)} \right ], \ \lambda >0.
 \end{equation}
As easily seen, one has (see e.g. \cite{abundo:MCAP2025}):
  \begin{equation} \label{QversusM}
\widehat Q_r(x, \lambda;x_R) = \frac 1 \lambda \left [ 1 - M _r (x, \lambda; x_R) \right ],
 \end{equation}
or
\begin{equation} \label{MversusQ}
M  _r (x, \lambda; x_R) = 1 - \lambda \widehat Q  _r(x, \lambda;x_R) .
 \end{equation}

By using \eqref{QversusM}, one also obtains  that $M _r (x, \lambda; x_R),$ as a function of $x ,$ satisfies the differential problem (cf. Eq. (2.3) of \cite{abundo:FPA2023}):
$$ {\cal L}M _r (x, \lambda; x_R) = \lambda M _r (x, \lambda; x_R),$$
namely
 \begin{equation} \label{eqLTtau}
 \frac {\sigma ^2 (x)} 2  \frac {\partial ^2} {\partial x^2 }M _r  (x, \lambda; x_R) + \mu (x) \frac {\partial} {\partial x } M  _r(x, \lambda; x_R)- (\lambda +r) M  _r(x, \lambda; x_R)
+ r M_r  (x_R, \lambda; x_R) = 0,
 \end{equation}
with boundary conditions
  \begin{equation}
M _r (0, \lambda; x_R) =1,
 \end{equation}
 in the case of  the FPT, and
 \begin{equation}
  M _r (0, \lambda; x_R) =M  (b, \lambda; x_R) =1,
   \end{equation}
 in the case of the FET. \bigskip

 If the LT of $\tau  (x,r)$
is finite for $\lambda$ belonging to a neighborhood of $0,$
then the $n-$th order moments of $\tau (x,r)$
exist  finite, and they are given by:
\begin{equation} \label{ODE1forM}
T_ n (x,r) := E \left [  \left ( \tau (x,r)
\right ) ^n \right ]  = (-1)^n \left [ \frac {\partial
^n } {\partial \lambda ^n } M_r (x, \lambda; x_R ) \right ] _ { \lambda =0} \ , n=1, 2 , \dots .
\end{equation}
The moments of $\tau (x,r)$ can be also obtained from  $\widehat Q _r(x, \lambda; x_R).$
For instance,  the first and second moments of $\tau (x,r)$ are  given by:
\begin{equation} \label{momentsoftauasQ}
E[\tau(x,r)] = \widehat Q_r(x, 0; x_R), \ \ E[\tau ^2(x,r)]= - 2 \frac \partial {\partial \lambda} \widehat Q _r(x, \lambda; x_R) |_ {\lambda =0}.
\end{equation}
By setting $T_0(x,r)=1$ and calculating  the $n-$th derivative with respect to $\lambda ,$ at $\lambda =0,$ of
both members  of   \eqref{eqLTtau}, we also obtain that
the $n-$th order moments $T_n(x,r)$  satisfy the ODEs (see e.g. \cite{abundo:TPMS2024} or \cite{abundo:FPA2023}):
\begin{equation}
{\cal L} T_n(x,r)= -n T_{n-1} (x,r) ,
\end{equation}
or
\begin{equation} \label{eqmoments}
 L T_n(x,r) = -n T_{n-1} (x,r) + r T_n(x,r) - r T_n(x_R,r) ,
\end{equation}
with  the boundary conditions  $T_n(0,r)=0$ in the case of the FPT, and $T_n (0,r)=T_ n(b,r)= 0$  for the FET (remember that the operator $L$ acts on $T_n(x,r),$ as a function of $x).$ \par\noindent
Note that for $r=0$, Eq. \eqref{eqmoments} becomes  the celebrated Darling and Siegert's equation  (\cite{darling:ams53}) for the moments of the first-passage time of a diffusion without resetting.
\par\noindent
In particular, $T_1(x,r)= E[\tau (x,r)]= \widehat Q_r (x,0; x_R),$
as a function of $x,$ is the solution of the differential problem with boundary conditions:
\begin{equation} \label{eqmeanFPT}
\begin{cases}
L T_1(x,r) - r  T_1 (x,r) + r T_1(x_R,r) = -1 ,\\
T_1(0,r) =0 \ ({\rm for \ the \ FPT}),  \\
T_1(0,r) = T_1(b,r) =0 \ ({\rm for \ the \ FET}),
\end{cases}
\end{equation}
whilst $E[\tau  ^2(x,r)]$ is the solution of the differential problem:
\begin{equation} \label{eqsecondFPT}
\begin{cases}
LT_2(x,r) - r  T_2(x,r) + r T_2(x_R,r) = -2 E[\tau (x,r)],  \\
T_2(0,r) =0 \ ({\rm for \ the \ FPT}), \\
T_2(0,r) = T_2(b,r) =0 \ ({\rm for \ the \ FET}).
\end{cases}
\end{equation}
A more convenient way to obtain the first and second moments of the  FPT  and the FET, which is alternative to solve the differential problems \eqref{eqmeanFPT} and \eqref{eqsecondFPT}, is to use the following:
\begin{Proposition} \label{LTandEFPT}
The Laplace transform of $\tau (x,r),$ i.e.  $M _r (x, \lambda; x_R) = E[e^ {-\lambda \tau (x,r)}],$  is given by:
\begin{equation} \label{MversusM0}
M_r (x, \lambda; x_R) = \frac {r M_0(x_R, \lambda +r) + \lambda M_0(x, \lambda +r )}  { \lambda +r M_0 (x_R, \lambda +r )}, \ \lambda >0,
  \end{equation}
where $M_0(x, \lambda) =E [e^ {- \lambda \tau (x,0)}]$ denotes the LT of the FPT (or FET) of the process without resetting $X(t)$ (that is, for $r=0).$ \par
Moreover, if the first and second moments of $\tau (x,r)$ exist finite, then:
\begin{equation} \label{expectedFPT}
E[\tau (x,r)]=  \frac {\widehat Q _0(x,r)} {1-r\widehat Q _0(x_R,r) }=
\frac 1 r  \left [\frac  {1-M_0(x,r)} {M_0(x_R,r)} \right ],
 \end{equation}
$$
E[\tau ^2 (x,r)]=  \frac 2 {r^2 M_0^2(x_R,r)} \Bigg \{ \frac {\partial M_0(x,r)} {\partial r}  \left [M_0(x_R,r) -1 \right ] +
$$
\begin{equation} \label{secondmomFPT}
-  \frac {\partial M_0(x_R,r)} {\partial r}  \left [ M_0(x,r) -1 \right ]
- \left [M_0(x_R,r)- r \frac {\partial M_0(x_R,r)} {\partial r}  -1 \right ]
\Bigg \},
 \end{equation}
 where $\widehat Q _0(x,r)$ is the LT of the survival function of the FPT (or FET) of $X(t).$ \bigskip

Furthermore, in terms of $\widehat Q_0,$ one has:
$$E[\tau ^2 (x,r)]=$$
\begin{equation} \label{secondmomFPT2}
= -  \frac 2  {\left (1-r \widehat Q_0 (x_R,r) \right )^2 }  \Bigg \{ \frac {\partial \widehat Q_0(x,r)} {\partial r}
-r \left [ \frac {\partial \widehat Q_0(x,r)} {\partial r} \widehat Q_0(x_R,r) -  \frac {\partial \widehat Q_0(x_R,r)} {\partial r} \widehat Q_0(x,r) \right ]  \Bigg \} .
\end{equation}
In particular, if $x=x_R, $ then:
 \begin{equation} \label{expectedFPTx=xR}
E[\tau (x,r)]= \frac 1 r \ \left [\frac  {1} {M_0(x,r)} -1 \right ] = \frac {\widehat Q _0(x,r)} {1-r \widehat Q _0(x,r) } .
\end{equation}
\begin{equation} \label{secondmomFPTx=xR}
E[\tau ^2(x,r)]= - \frac 2 {r^2 M_0^2(x,r)} \left (M_0(x,r)- r \frac {\partial M_0(x,r)} {\partial r}  -1 \right ).
\end{equation}
or
\begin{equation} \label{secondmomFPT2x=xR}
E[\tau ^2(x,r)]= -   \frac {2 } {\left (1-r \widehat Q_0 (x,r) \right )^2 } \cdot  \frac {\partial \widehat Q_0(x,r)} {\partial r}  .
\end{equation}
\end{Proposition}
{\it Proof.} From \eqref{renewalforQ2} and \eqref{MversusQ}, we soon obtain \eqref{MversusM0}; since
$E[\tau (x,r)]= \widehat Q  _r(x, \lambda;x_R),$ formula \eqref{expectedFPT}  follows
by using  \eqref{renewalforQ2}. \par\noindent
The expression for the second moment is obtained by making calculations in the formula \par\noindent
$E[\tau ^2(x,r)] = -2 \frac {\partial \widehat Q _r (x, \lambda; x_R)} {\partial \lambda} |_{\lambda =0},$ and using \eqref{renewalforQ2}, \eqref{QversusM}, \eqref{MversusQ}. \par\noindent
Finally, formulae \eqref{expectedFPTx=xR}, \eqref{secondmomFPTx=xR} and \eqref{secondmomFPT2x=xR} follow, by taking $x=x_R.$ \par
\hfill  $\Box$

\begin{Remark}
Proposition \ref{LTandEFPT} allows to write  the LT of the FPT (or FET) of the process with resetting ${\cal X}(t)$, in terms of the LT of the FPT (or FET)  of the underlying process without resetting
$X(t).$ Moreover, it provides the first and second moments of the  FPT (or FET) of ${\cal X}(t),$ in terms  of the LT of  the FPT (or FET) of $X(t). $
In some cases, this approach is more convenient, from a computational point of view; in fact, the solutions of the differential problems \eqref{eqmeanFPT}, \eqref{eqsecondFPT} for the first two moments of
$\tau(x,r)$ were obtained, in the case of drifted BM with resetting, in \cite{abundo:MCAP2025}, \cite{abundo:FPA2023} in a rather laborious way,
\end{Remark}

\section{Some examples of explicit calculation of the expected FPT}
In this section, we provide explicit calculations of the expected FPT through zero, for some
diffusions with resetting ${\cal X}(t).$
\subsection{Drifted Brownian motion with resetting}
Let $X(t) = x+  \eta t + B_t$ be BM with constant drift $\eta ,$ starting from $x >0,$ and  ${\cal X}(t)$ the corresponding process with resetting.
For non-negative reset position $x_R,$ we will calculate the mean of the  FPT  of ${\cal X}(t)$ through zero; in fact, it is finite for $r \neq 0$ (see e.g. \cite{abundo:FPA2023}).   \par\noindent
We recall that the LT of the FPT of $X(t)$ through zero is (see e.g. \cite{borodin:1996}):
\begin{equation} \label{LTFPTdriftBM}
M_0(x, \lambda) = E[e^ {- \lambda \tau (x,0)}] = e^ {-x \left (\eta +  \sqrt {\eta ^2 + 2 \lambda}  \ \right )}, \ \lambda >0 .
 \end{equation}
From \eqref{LTFPTdriftBM} and \eqref{MversusM0} of Proposition \ref{LTandEFPT}, we get
\begin{equation} \label{LTFPTdriftresetBM}
M_r(x, \lambda; x_R) = E[e^ {- \lambda \tau (x,r)}] = \frac {\lambda e^ {-x \left (\eta +  \sqrt {\eta ^2 + 2 (\lambda +r)} \ \right )} +r e^ {-x_R \left ( \eta + \sqrt {\eta ^2 + 2 (\lambda +r)} \ \right )}}
{\lambda +r e^ {-x_R \left ( \eta + \sqrt {\eta ^2 + 2 ( \lambda +r)} \ \right )} }, \ \lambda >0 .
 \end{equation}
 Moreover, by using Eq. \eqref{expectedFPT} and \eqref{secondmomFPT}, we obtain  the first and second moments of the  FPT of drifted BM with resetting. \par\noindent
 For instance, for $r >0$ the expected FPT is:
\begin{equation} \label{EFPTBMreset}
E[\tau (x,r)] =  \frac 1 r e^ {x_R \left (\eta + \sqrt {\eta ^2 + 2 r} \ \right )} \left ( 1- e^ {-x \left (\eta + \sqrt {\eta ^2 + 2 r} \  \right )} \right ).
\end{equation}
The expression of $E[\tau ^2(x,r)]$ is rather long and we omit to write it, here. \par\noindent
These formulae were obtained in \cite{abundo:FPA2023} (see Eqs. (4.1), (4.3) and (4.4) therein) in a rather laborious way, by solving the differential problems
\eqref{eqLTtau}, \eqref{eqmeanFPT} and \eqref{eqsecondFPT}
with $Lu(x) = \frac 1 2  u''(x) +  \eta u'(x).$

\subsection{Ornstein-Uhlenbeck process with resetting} \par\noindent
Let $X(t)$ be Ornstein-Uhlenbeck (OU) process, starting from $x >0,$ namely the solution of the SDE:
\begin{equation} \label{OUSDE}
dX(t) = - \mu X(t) dt + \sigma dB_t,
 \end{equation}
 with initial condition $X(0) = x > 0,$
where $\mu$ and $\sigma$ are positive constants.  While for BM the drift coefficient is a constant function, now it is proportional to $X(t)$ through the negative constant $- \mu.$\par\noindent
As well-known, it holds the representation
\begin{equation}
X(t)= e^ {- \mu t} \left ( x + B(\rho (t))  \right ),
\end{equation}
where
$
\rho (t)= \frac {\sigma ^2} {2 \mu} \left (e^ {2 \mu t} -1 \right ).
$
\par\noindent
Let ${\cal X}(t)$ be the corresponding OU process with resetting.
For non-negative reset position $x_R,$ we will calculate the mean of the  FPT  of ${\cal X}(t)$ through zero.   \par\noindent
The LT $M_0(x, \lambda) = E[e^ {- \lambda \tau (x,0)}]$ of the FPT of OU process $X(t)$ through zero, under the condition that $X(0)=x >0 ,$ satisfies the differential equation \eqref{eqLTtau} with $r=0$ and $Lu(x) = \frac {\sigma ^2} 2  u''(x) - \mu x u'(x),$
namely:
\begin{equation} \label{diffeqLTOU}
\frac {\sigma ^2 } 2  \frac {d^2 u} {dx^2} - \mu x \frac {d u} {dx} - \lambda u =0 ,
\end{equation}
with boundary conditions $u(0)=1$ and $\lim _{x \rightarrow + \infty} u(x)=0.$ \par\noindent
Its explicit solution is (see e.g. \cite{capric:71} or \cite{rino:99} ):
\begin{equation} \label{LTFPTOU}
M_0(x, \lambda) =   e^  {\frac {\mu x^2 } {2 \sigma ^2} } \  \frac {D_ {- \lambda / \mu} \left (x \sqrt {2 \mu / \sigma ^2} \right )} {D_ {- \lambda / \mu} ( 0)}, \ \lambda >0 ,
\end{equation}
where $D _ \nu $ is the Parabolic Cylinder function (see \cite{grad:80}), which satisfies  the ODE:
\begin{equation}
\frac { d^2 D_\nu (z)} {dz^2 } + \left ( \nu + \frac 1 2 - \frac { z^2} 4 \right ) D_ \nu (z)=0, \
D_ \nu (0) =
\frac { 2^ {\nu /2} \sqrt \pi } { \Gamma ( \frac { 1 - \nu} {2 } ) } .
\end{equation}
Actually, formula \eqref{LTFPTOU} is obtained from the analogous ones in \cite{capric:71} or \cite{rino:99}, by taking $-x$ in place of $x.$ In fact, formulae in \cite{capric:71} or \cite{rino:99} concern the FPT of the OU process through a boundary $S,$ when starting from $x <S;$ thanks to the symmetry of OU, the FPT through zero, when starting from $x >0$ is nothing but the FPT through zero, when starting from $-x.$ \par\noindent
By using \eqref{LTFPTOU} and  Proposition \ref{LTandEFPT}, one can  get the expression of the LT of the FPT of ${\cal X}(t);$
from \eqref{expectedFPT}, we easily obtain that the expected FPT through zero of OU with resetting is:
 \begin{equation} \label{EFPTOUreset}
 E[\tau (x,r)] = \frac 1 r \left [ \left (1- e^  {\frac {\mu x^2 } {2 \sigma ^2} } \  \frac {D_ {- r / \mu} \left (x \sqrt {2 \mu / \sigma ^2} \right )} {D_ {- r / \mu} ( 0)}  \right ) \cdot
 \frac 1 {e^  {\frac {\mu x^2 } {2 \sigma ^2} } \  \frac {D_ {- r / \mu} \left (x_R \sqrt {2 \mu / \sigma ^2} \right )} {D_ {- r / \mu} ( 0)}    }    \right ].
 \end{equation}
 In particular, if $x=x_R:$
  \begin{equation} \label{EFPTOUx=XRreset}
  E[\tau (x,r)] = \frac 1 r \ \left [ \frac {D_ {- r / \mu} ( 0)}  {D_ {- r / \mu} \left (x \sqrt {2 \mu / \sigma ^2} \right )}  \ e^  {\frac {- \mu x^2 } {2 \sigma ^2} } -1 \right ].
  \end{equation}
\par
Note that, for $r=0$  an alternative integral formula holds for $ E[\tau (x,0)],$ that does not involve special functions (see e.g. \cite{abundo:OUarea}):
 \begin{equation} \label{altermeantau}
  E[\tau (x,0)]=
\frac 1 \mu \int _0 ^ { + \infty } e^ { - \frac {\sigma ^2 } {2 \mu } \ \frac {y^2} 2 } \ \frac {1 - e^ {-xy}}  y \ dy .
\end{equation}
From \eqref{secondmomFPT} and \eqref{LTFPTOU}, one could also obtain the second moment $E[\tau ^2(x,r)],$ by  calculations.
\subsection{CIR process with resetting} \par\noindent
For $\mu >0$ and $\sigma >0,$ let  $X(t)$ be the solution of the SDE:
\begin{equation} \label{CIReq}
dX(t)= (\sigma ^2 -2 \mu X(t)) dt + 2 \sigma \sqrt {X(t)} dB_t, \ X(0) = x >0.
\end{equation}
The diffusion $X(t)$ turns out to be non-negative for all $t \ge 0$ and describes a CIR model
of the type
$$ dX(t)= ( \beta - 2 \alpha X(t)) dt + \delta \sqrt {X(t)} dB_t ,$$
for which the Feller condition for accessibility of the zero boundary, i.e. $\frac {2 \beta} {\delta ^2} <1,$
is satisfied, since $2 \sigma ^2 < 4 \sigma ^2.$ \par\noindent
As easily seen,  $Y(t):= \sqrt {X(t)}$ is OU process, driven by the SDE \eqref{OUSDE} with $Y(t)$ in place of $X(t).$ Thus, the FPT, $\tau (x,0),$  of $X(t)$ through zero, under the condition that $X(0)= x >0,$ is nothing but
the FPT of $Y(t)$ through zero, under the condition that $Y(0)= \sqrt x ,$ and therefore
\begin{equation} \label{LTFPTCIR}
M_0(x, \lambda)= E[e^ {-\lambda \tau (x, 0)}] = M_0 ^ {OU} (\sqrt x, \lambda), \ \lambda >0,
\end{equation}
where $M_0 ^ {OU} ( z, \lambda)$ is the LT of the FPT of OU starting from $z,$  which is given by  \eqref{LTFPTOU}. \par\noindent
Then, if ${\cal X}(t)$ is the corresponding CIR process with resetting,
by the previous result it follows that the expected FPT through zero, i.e. $E[\tau (x,r)],$  is obtained from
\eqref{EFPTOUreset}, with $\sqrt x$ in place of $x$ and $\sqrt {x_R}$ in place of $x_R,$ namely:
 \begin{equation} \label{EFPTCIR}
 E[\tau (x,r)] = \frac 1 r \left [ \left (1- e^  {\frac {\mu x } {2 \sigma ^2} } \  \frac {D_ {- r / \mu} \left (\sqrt {2 \mu x / \sigma ^2} \right )} {D_ {- r / \mu} ( 0)}  \right ) \cdot
 \frac 1 {e^  {\frac {\mu x } {2 \sigma ^2} } \  \frac {D_ {- r / \mu} \left (\sqrt {2 \mu x_R / \sigma ^2} \right )} {D_ {- r / \mu} ( 0)}    }    \right ].
 \end{equation}
As in the case of OU with resetting, from \eqref{secondmomFPT} and \eqref{LTFPTCIR}, one could also obtain the second moment $E[\tau ^2(x,r)],$ by  calculations.
\subsection{Diffusions with resetting conjugated to BM} \par\noindent
A diffusion $X(t),$ starting from $X(0)=x,$ is said to be {\it conjugated} to BM, if
there exists an increasing differentiable function $v(x)$ with $v(0) = 0,$ such
that $X(t)= v^{-1} \left (B_t + v( x ) \right ),$ for any $t \ge 0.$ (see \cite{abundo:stapro12}).
Two simple examples of diffusions conjugated to BM are the following: \par
${\bf \bullet}$  the  Feller process $X(t),$ driven by the SDE $dX(t)= \frac 1 4 dt + \sqrt { X(t) } \ dB_t  , \ X(0)= x >0,$  which is conjugated to
BM via the the function $v(x)= 2 \sqrt x ;$
\par
${\bf \bullet}$ the Wright-Fisher \ like \ process, driven by the SDE \par\noindent
$dX(t)=
(\frac 1 4 - \frac 1 2 X(t))dt + \sqrt { X(t) (1- X(t)) } \ dB_t  , \ X(0)= x  \in (0,1) ,$ which is conjugated to
BM via the the function $v(x)= 2 \arcsin \sqrt x.$ \par
If $X(t)$  is conjugated to BM via the function $v$, then
the FPT  of $X(t)$ through zero, under the condition that $X(0)= x,$ is nothing but
the FPT of BM through zero, when starting from $v(x),$  and therefore the LT of $\tau (x,0)$ is obtained from  \eqref{LTFPTdriftBM}, with $\eta =0$ and $v(x)$ in place of $x.$
\par
If ${\cal X}(t)$ is the corresponding diffusion with resetting,  then $v({\cal X}(t))$ is BM with resetting, so  the expected FPT  through zero is given by
\eqref{EFPTBMreset}, with $\eta =0, \ v(x)$ in place of $x,$ and  $v(x_R)$ in place of $x_R,$ namely:
\begin{equation} \label{EFPTconjBMreset}
E[\tau (x,r)] =  \frac 1 r e^ {v(x_R) \sqrt { 2 r} } \left ( 1- e^ {-v(x)  \sqrt {2 r} } \right ).
\end{equation}
As in the case of BM with resetting, from \eqref{secondmomFPT} and \eqref{LTFPTdriftBM} (with $\eta =0$ and $v(x)$ in place of $x),$ one could also obtain the second moment $E[\tau ^2(x,r)].$

\section{Some examples of explicit calculation of the expected FET}
In this section, we present some examples of explicit calculations of the expected FET of a diffusion with resetting ${\cal X}(t) $ from the interval $(0,b).$
\subsection{Drifted Brownian motion with resetting}
For  drifted BM $X(t)= x + \eta t + B_t ,$  let us consider the FET
$\tau (x,0) = \inf \{t >0: X(t) \notin (0,b) \};$ then, the survival function of $\tau (x,0),$ i.e. $Q _0(x,t)= P[ \tau  (x,0)>t],$ satisfies  Eq.
\eqref{eqLTsurvivalprob} with $r=0$ and $L(u)= \frac 1 2 u''(x) + \eta u'(x);$ its explicit solution is (see e.g. \cite{abundo:MCAP2025}):
\begin{equation} \label{explicit LTQdriftBMreset}
\widehat Q _0(x, \lambda) = \frac {e^{-b \eta}\sinh(b \beta _ \lambda) - e^ {- (b+x) \eta }\sinh(\beta _ \lambda (b-x)) - e^ {- \eta x}\sinh (\beta _ \lambda x) }
{\lambda e^{- b \eta} \sinh (\beta _ \lambda b )} ,
\end{equation}
where $\beta _ \lambda = \sqrt {\eta ^2 + 2\lambda}$
(we omit to explicit the dependence on $\eta , $ for simplicity). \par
Moreover, the LT of $\tau  (x,0)$ satisfies Eq. \eqref{eqLTtau} with $r=0$ and the same infinitesimal generator $L(u);$ its explicit solution is  (see e.g. \cite{borodin:1996}, Eq. 3.01, pg. 233):
\begin{equation} \label{LTtaudriftBMreset}
M_0 (x, \lambda) = 1- \lambda \widehat Q _0(x, \lambda)=
\frac {e^ {- \eta x}\sinh (\beta _ \lambda (b-x)) + e^ { (b-x) \eta }\sinh(\beta _ \lambda x) } {\sinh(b \beta _ \lambda)} ,
\end{equation}
By using \eqref{expectedFPT} and  \eqref{explicit LTQdriftBMreset}, we obtain that the expected FET of the corresponding drifted BM with resetting is:
$$ E[\tau (x,r)] = $$
\begin{equation} \label{EtauBMdriftreset}
= \frac {e^ {-b \eta }\sinh \left (b \sqrt {\eta ^2 +2r} \ \right )- e^ {- \eta x}\sinh \left (x \sqrt {\eta ^2 +2r} \ \right )
- e^ {-(b+x) \eta}\sinh \left ((b -x) \sqrt {\eta ^2 +2r} \ \right )   }  {r \left [ e^ {-x_R \eta}\sinh \left (x_R \sqrt {\eta ^2 +2r} \ \right ) + e^ {-(b+x_R) \eta}\sinh \left ((b - x_R) \sqrt {\eta ^2 +2r}  \ \right ) \right ]}.
\end{equation}
This formula was already found in \cite{abundo:MCAP2025} (Eq. (2.42)), by solving the differential problem \eqref{eqmeanFPT}
with $Lu= \frac 1 2 u'' + \eta u',$ though its derivation was rather complicated. \par\noindent
As in the case of the FPT, from \eqref{secondmomFPT} and \eqref{LTtaudriftBMreset}, one could also obtain the second moment of the FET of drifted BM with resetting.

\subsection{OU with resetting} \par\noindent
The LT of the FET of OU process from the interval $(0, b)$ satisfies Eq.
\eqref{diffeqLTOU}, with the boundary conditions $u(0, \lambda )=1$ and  $u(b, \lambda )=1 .$ The general solution of this ODE is a linear combination of two independent solutions, which are:
\begin{equation} \label{functionsui}
u_1(x, \lambda)= e^  {\frac {\mu x^2 } {2 \sigma ^2} } D_ {- \lambda / \mu} \left (-x \sqrt {2 \mu / \sigma ^2} \right ), \  u_2(x, \lambda )= e^  {\frac {\mu x^2 } {2 \sigma ^2} } D_ {- \lambda / \mu} \left (x \sqrt {2 \mu / \sigma ^2} \right ),
\end{equation}
where  $D _ \nu $ is the Parabolic Cylinder function  (see Section 3.2).
Thus, the LT of the FET of OU process can be written as
\begin{equation} \label{LTFETOU}
M_0(x, \lambda) = E \left [e^ {- \lambda \tau (x,0)} \right ] = c_1 u_1(x, \lambda) + c_2 u_2 (x, \lambda),
\end{equation}
where $c_1$ and $ c_2$ are constants (with respect to $x)$ to be found. By imposing the boundary conditions $M_0(0, \lambda) = M_0(b, \lambda) =1,$ one finds
\begin{equation} \label{c1}
c_1= \frac {e^  {\frac {\mu b^2 } {2 \sigma ^2} } D_ {- \lambda / \mu} \left (b \sqrt {2 \mu / \sigma ^2} \right ) -
D_ {- \lambda / \mu} \left (0 \right )} {D_ {- \lambda / \mu} \left ( 0 \right ) e^  {\frac {\mu b^2 } {2 \sigma ^2} } D_ {- \lambda / \mu} \left (b \sqrt {2 \mu / \sigma ^2} \right ) -
e^  {\frac {\mu b^2 } {2 \sigma ^2} } D_ {- \lambda / \mu} \left (-b \sqrt {2 \mu / \sigma ^2} \right ) D_ {- \lambda / \mu} \left ( 0 \right ) },
\end{equation}
\begin{equation}  \label{c2}
c_2= \frac { D_ {- \lambda / \mu} \left ( 0 \right ) - e^  {\frac {\mu b^2 } {2 \sigma ^2} } D_ {- \lambda / \mu} \left (-b \sqrt {2 \mu / \sigma ^2} \right ) }   {D_ {- \lambda / \mu} \left ( 0 \right ) e^  {\frac {\mu b^2 } {2 \sigma ^2} } D_ {- \lambda / \mu} \left (b \sqrt {2 \mu / \sigma ^2} \right ) -
e^  {\frac {\mu b^2 } {2 \sigma ^2} } D_ {- \lambda / \mu} \left (-b \sqrt {2 \mu / \sigma ^2} \right ) D_ {- \lambda / \mu} \left ( 0 \right ) }.
\end{equation}
Thus, from \eqref{LTFETOU}, the explicit form of the LT of the FET  $\tau (x,0)$ of OU process follows. \par
Note that formula \eqref{LTFPTOU}, for the LT of the FPT of OU process (without resetting) through zero, is also obtained from
\eqref{LTFETOU} by imposing the boundary conditions:
$$M_0(0, \lambda)=1, \  {\rm and}  \ \lim _ {x \rightarrow + \infty } M_0(x, \lambda)=0,$$
which provide $c_1=0$ and $c_2= 1/ D_{- \lambda / \mu} (0).$
\par
Turning back to consider the FET, by using \eqref{expectedFPT} we get that the expected FET for OU with resetting is:
\begin{equation} \label{meanFETOU}
E[\tau (x,r)] = \frac 1 r \left [ \frac {1- c_1u_1(x,r) -c_2 u_2(x,r)} {c_1u_1(x_R,r) +c_2 u_2(x_R,r)}\right ],
\end{equation}
where $u_i (x,r)$ and $c_i$ are given by \eqref{functionsui}, and  \eqref{c1}, \eqref{c2}.
From \eqref{secondmomFPT} one could also obtain  $E[\tau ^2(x,r)].$
However, both  expressions of first and second moments of the FET of OU with resetting are a bit more complicated, than the corresponding formulae, concerning  the FPT (see e.g. \eqref{EFPTOUreset} for the first moment).

\subsection{CIR process with resetting} \par\noindent
For $\mu >0$ and $\sigma >0,$ let  $X(t)$ be the solution of the SDE \eqref{CIReq};
then, the diffusion $Y(t)= \sqrt {X(t)}$ is OU process.
So, the expected FET of CIR process with resetting from $(0,b)$ is easily obtained by using \eqref{meanFETOU}
 with $\sqrt x$ in place of $x$ and $\sqrt {x_R}$ in place of $x_R,$ namely:
\begin{equation} \label{meanFETCIR}
E[\tau (x,r)] = \frac 1 r \left [ \frac {1- c_1u_1(\sqrt x,r) -c_2 u_2(\sqrt x,r)} {c_1u_1(\sqrt {x_R},r) +c_2 u_2(\sqrt {x_R},r)}\right ],
\end{equation}
where $u_i(x,r)$ are as in \eqref{functionsui}, but $c_i$ are obtained from \eqref{c1}, \eqref{c2}, by replacing $b$ with $\sqrt b.$ \par\noindent
As in the case of the FPT, by using \eqref{secondmomFPT} one could also obtain  the second moment of the FET of CIR process with resetting.

\subsection{Diffusions with resetting conjugated to BM}
Let   $X(t)$  be conjugated to BM via the function $v$, then
the FET $\tau (x,0)$  of $X(t)$ from the interval $(0,b),$ under the condition that $X(0)= x,$ is nothing but
the FET of BM from the interval $(0,v(b)),$  when starting from $v(x),$  and therefore the LT of $\tau (x,0)$ is obtained from  \eqref{LTtaudriftBMreset}, with $\eta =0 , \ v(x)$ in place of $x,$
and $v(b) $ in place of $b.$
\par\noindent
If ${\cal X}(t)$ is the corresponding diffusion with resetting,  then $v({\cal X}(t))$ is BM with resetting, so by the result concerning BM  the expected FET from $(0,b), \  E[\tau (x,r)],$ is given by
\eqref{EtauBMdriftreset}, with $\eta =0, \ v(x)$ in place of $x,$ \ $v(b) $ in place of $b,$ and  $v(x_R)$ in place of $x_R,$ namely:

\begin{equation} \label{EtauBMconjreset}
E[\tau (x,r)] = \frac {\sinh \left (v(b) \sqrt {2r} \ \right )- \sinh \left (v(x) \sqrt {2r} \ \right )
- \sinh \left ((v(b) - v(x)) \sqrt {2r} \ \right )   }  {r \left [ \sinh \left (v(x_R) \sqrt {2r} \ \right ) + \sinh \left ((v(b) - v(x_R)) \sqrt {2r}  \ \right ) \right ]}.
\end{equation}
As in the case of the FPT, by using \eqref{secondmomFPT} one could also obtain  the second moment of the FET of  ${\cal X}(t).$

\section{Minimization of the expected FPT and the FET}
In this section we will deal with the minimization of the expected FPT and  the expected FET of a diffusion with resetting, with respect to the resetting rate $r;$ actually,
we will consider only  drifted BM and OU, with resetting. For drifted BM with resetting we will treat both the cases of the FPT and the FET, while for OU with resetting we
will limit ourselves to the FPT,
because  formula \eqref{meanFETOU} for the expected FET is rather computationally complicated. We will not report explicitly the
minimization of the expected FPT and FET for CIR and  processes
conjugated to BM, because they can be  traced back to those of OU and BM.

\subsection{The case of drifted Brownian motion with resetting}
\subsubsection{Optimization of the FPT}
Let $X(t) = x+  \eta t + B_t$ be BM with drift $\eta,$ starting from $x >0,$ and  ${\cal X}(t)$ the corresponding process with resetting;
we suppose that the reset position $x_R$ is non-negative.
The expectation of  the FPT,  $\tau (x,r) = \inf \{t>0: {\cal X} (t) = 0 \ | \ {\cal X} (0)=X(0) =x  \},$  is given by  \eqref{EFPTBMreset}. \bigskip

\noindent {\bf (a) BM with resetting $(\eta =0)$} \par
One has:
\begin{equation} \label{Etau}
T(x,r):= E[ \tau (x,r)] =
\begin{cases}
\frac 1 r e^ {x_R \sqrt {2r}} \left ( 1- e^{-x \sqrt {2r}} \right ), & x,r >0 \\
0, &  x=0, \ r \ge 0 \\
+ \infty , &  x >0, \ r=0.
\end{cases}
\end{equation}
\par
Unlike the case of BM without resetting $(r=0),$ the expectation of the FPT, $T(x,r) ,$ results to be finite for all $x >0$ and $r>0$ (see e.g. \cite{abundo:FPA2023}), and:
\begin{equation}
\lim _ {r \rightarrow 0^+} T(x,r) = + \infty, \ {\rm as \ well \ as } \ \lim _ {r \rightarrow + \infty } T( x,r) = + \infty.
\end{equation}
As easily seen,  for fixed reset position $x_R >0$ and starting point $x >0,$ the expected FPT,  as a function of $r,$ attains the unique global minimum at a value
$
r_m(x)= arg \left ( \min _ {r \ge 0 } T (x,r) \right ).
$
Our goal is to find $r_{m}(x)$ and the
minimum expected FPT, $m(x)= T ( x,r_{m}(x)) ,$ for fixed $x_R$ and $x >0.$  In this way, $T(x,r)$ can be reduced, by resetting the process at the optimal value  $r_m (x).$\bigskip

\noindent {\bf (i) The case when $x=x_R$ } \par\noindent
For $x = x_R >0$ and $r >0 ,$ formula \eqref{Etau} becomes:
\begin{equation}
T(x,r)= \frac 1 r  \left ( e^{x \sqrt {2r}} -1 \right ).
\end{equation}
For small $x >0,$ one gets $T(x,r) \simeq \frac {x \sqrt 2} {\sqrt r}$ which is decreasing, as a function of $r,$ moreover
$T(x,r) \rightarrow + \infty,$ as $r \rightarrow 0^+,$  while $T(x,r) \rightarrow 0,$ as $r \rightarrow + \infty .$  \par\noindent
For large $x >0,$ it results $T(x,r) \simeq \frac {e^ {x \sqrt {2r}}} r $ which attains its global minimum at $r=2/x^2 .$ Therefore, for large, but finite $x,$ the choice of
the reset rate $r=2/x^2$ expedites the FPT. \par
In general, since the equation $\frac \partial {\partial r} T( x,r) =0$ cannot be solved analytically, in order to find the value $r_{m}(x)$ at which the minimum of $T( x,r)$ is attained,
we have to solve numerically it. \par\noindent
In Table \ref{tab1} below, we report the values of $r_{m}(x)$ and $m(x) = \min _ {r \ge 0 } T (x,r)= T(x, r_{m}(x)),$ numerically obtained  by
the secant method, for some values of $x=x_R >0.$ It results approximately
$r_m(x)= 1.27 / x^{2}$ and $m(x)= 3.088 \ x^2$ (see also \cite{Pal:19}).
We see that, as $x$ increases, $r_m(x)$ decreases and  it approximates zero, for large $x,$  whereas $m(x)$ increases. For small $x,$ a large value of the reset rate $r$
is needed to minimize the expected FPT, and the corresponding minimum of the expected FPT is small,
 while, for large $x$ a small value of $r$ is required, but the corresponding minimum of the expected FPT turns out to be large.\par\noindent
 \begin{table}[!h]
\begin{center}
\begin{tabular}{cccc}
\hline
 $x$ & $r_{m}(x)$ & $m(x)$  \\
\hline
0.1    & 126.980 & 0.030 \\
0.5    & 5.079 & 0.772 \\
1.     & 1.269 & 3.088 \\
2.     & 0.317 & 12.353 \\
3.     & 0.141 & 27.79 \\
5.     & 0.050 & 77.206 \\
10.     & 0.012 & 308.827 \\
\hline
\end{tabular}
\end{center}
\caption{Minimization of the expected FPT, $T(x,r)$, of BM with resetting: the table report the values of $r_{m}(x)$ and $m(x)=T(x,r_{m}(x))$ numerically obtained, for some values of $x =x_R >0.$  }
\label{tab1}
\end{table}
\bigskip

\noindent {\bf (ii) The case when $x \neq x_R $} \par\noindent
In this case, we observe a more intricate scenario.
Even now we have to numerically calculate the value $r_{m}(x)$ at which the minimum of $T( x,r)$ is attained $(T(x,r)$ is now given by \eqref{Etau}). \par\noindent
As an example, in  Table \ref{tab2}, we report the values of $r_{m}(x)$ and $m(x) = \min _ {r \ge 0 } T (x,r)= T(x, r_{m}(x))$, numerically obtained
by the secant method, for several values of $x >0,$ with fixed $x_R =1.$  Note that, as it must be, for $x=1$ the values of $r_m(x)$ and $m(x)$ coincide with those of Table \ref{tab1}. \par\noindent

\begin{table}[!h]
\begin{center}
\begin{tabular}{cccc}
\hline
 $x$ & $r_{m}(x)$ & $m(x)$  \\
\hline
0.0001 & 0.5000 & 0.00054 \\
0.001  & 0.5005 & 0.05430 \\
0.01   & 0.5050 & 0.05409 \\
0.1    & 0.5529 & 0.51671 \\
0.3    & 0.6780 & 1.39345 \\
0.5    & 0.8289 & 2.07535 \\
0.9    & 1.1812 & 2.94998 \\
1.     & 1.2698 & 3.08827 \\
1.5    & 1.6323 & 3.48334 \\
2.     & 1.7323 & 3.57000 \\
2.5    & 1.8000 & 3.65000 \\
3.     & 1.9691 & 3.68515 \\
5.     & 1.9990 & 3.69436 \\
7.     & 1.9990 & 3.69505 \\
\hline
\end{tabular}
\end{center}
\caption{Minimization of the expected FPT, $T(x,r)$, of BM with resetting: the table report the values of $r_{m}(x)$ and $m(x)=T(x,r_{m}(x))$ numerically obtained  for some values of $x >0,$ with fixed $x_R =1.$
The functions $r_m(x)$ and $m(x)$ are both increasing.}
\label{tab2}
\end{table}

For fixed $x_R >0,$ and $x>0$ close to zero,
$ T(x,r) $ behaves as
$\frac {x \sqrt 2} {\sqrt {r}}e^ {x_R \sqrt {2r}}.$
By calculating the derivative of this function with respect to $r$ and imposing it to be zero, one obtains
$r_{m}(x) = 1/2 x_R^2 .$ Instead, for large $x >0,$ \ $T(x,r) $ behaves as $e^ {x_R \sqrt {2r}}/r,$ which attains its global minimum at $r_{m}=2/x_R^2 .$ \par\noindent
Actually, for fixed $x_R$ the argument  $r_{m}(x)$ of  $\min _{r \ge 0} T(x,r) $ turns out to be  an increasing function of $x >0,$ and
$r_{m}(x)$  takes values within a range $(\alpha, \beta),$ where $\alpha = \lim _ {x \rightarrow 0 ^+} r_{m}(x) = \frac 1 {2 x_R^2}$ and $\beta = \lim _ {x \rightarrow + \infty} r_{m}(x)= \sup _{x >0} r_{m}(x) = \frac 2 {x_R^2}= 4 \alpha $
(therefore, $\alpha$ and $\beta$ turn out to be decreasing functions of the reset position $x_R).$ \par\noindent
Since $T(0,r)$ is  zero for every $r \ge 0,$ we can set $r_m(0)=0;$ instead, $\lim _{x \rightarrow 0^+} r_m(x)= \alpha >0,$ hence the function $r_m(x)$ has a jump discontinuity point at
$x=0.$
\par\noindent
Meanwhile, for $x \ge 0$ the minimum $m(x)= \min _ {r \ge 0} T(x,r) $ also increases from $0$ to $e^ {x_R \sqrt {2 \beta}}/ \beta .$  \par
\begin{figure}[!h]
\centering
\includegraphics[height=0.35 \textheight]{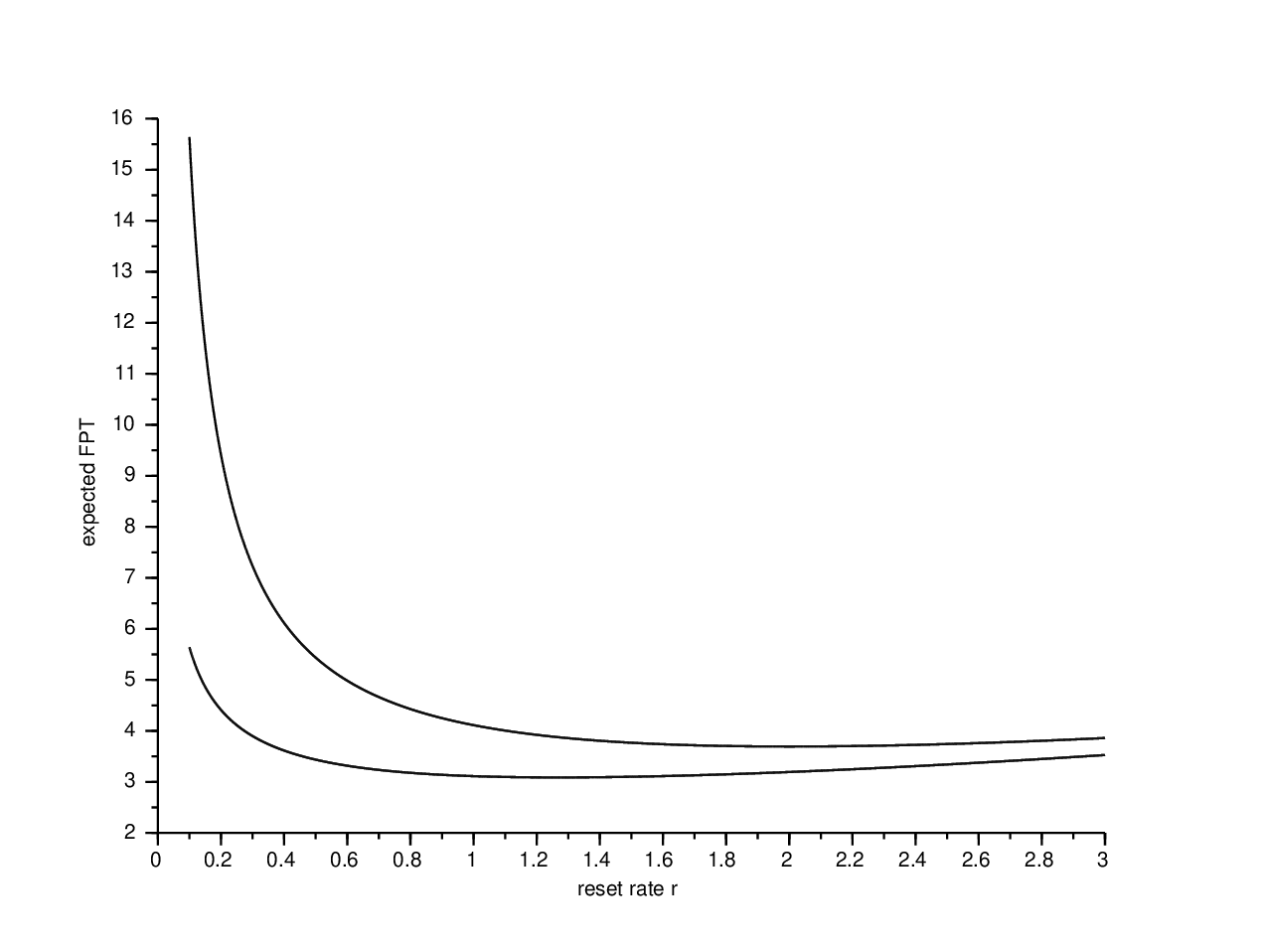}
\caption{For BM with resettings, the figure shows the graphs of the expected FPT $T(x,r),$  as  functions of $r>0,$ for fixed $x_R=1$ and  $x=1$ (lower curve),   $x=200$ (higher curve)
(on the horizontal axes $r).$
}
\label{FPTminbis}
\end{figure}

For fixed $x_R=1,$ Figure \ref{FPTminbis} shows two examples of the graphs of the expected FPT $T(x,r) ,$ as functions of $r>0,$ for  $x=1$ and $x=200;$ we see that the first curve attains the minimum at $r \simeq 1.269,$ while the second one attains the minimum at $r \simeq 2.$ \par
\begin{figure}[!h]
\centering
\includegraphics[height=0.35 \textheight]{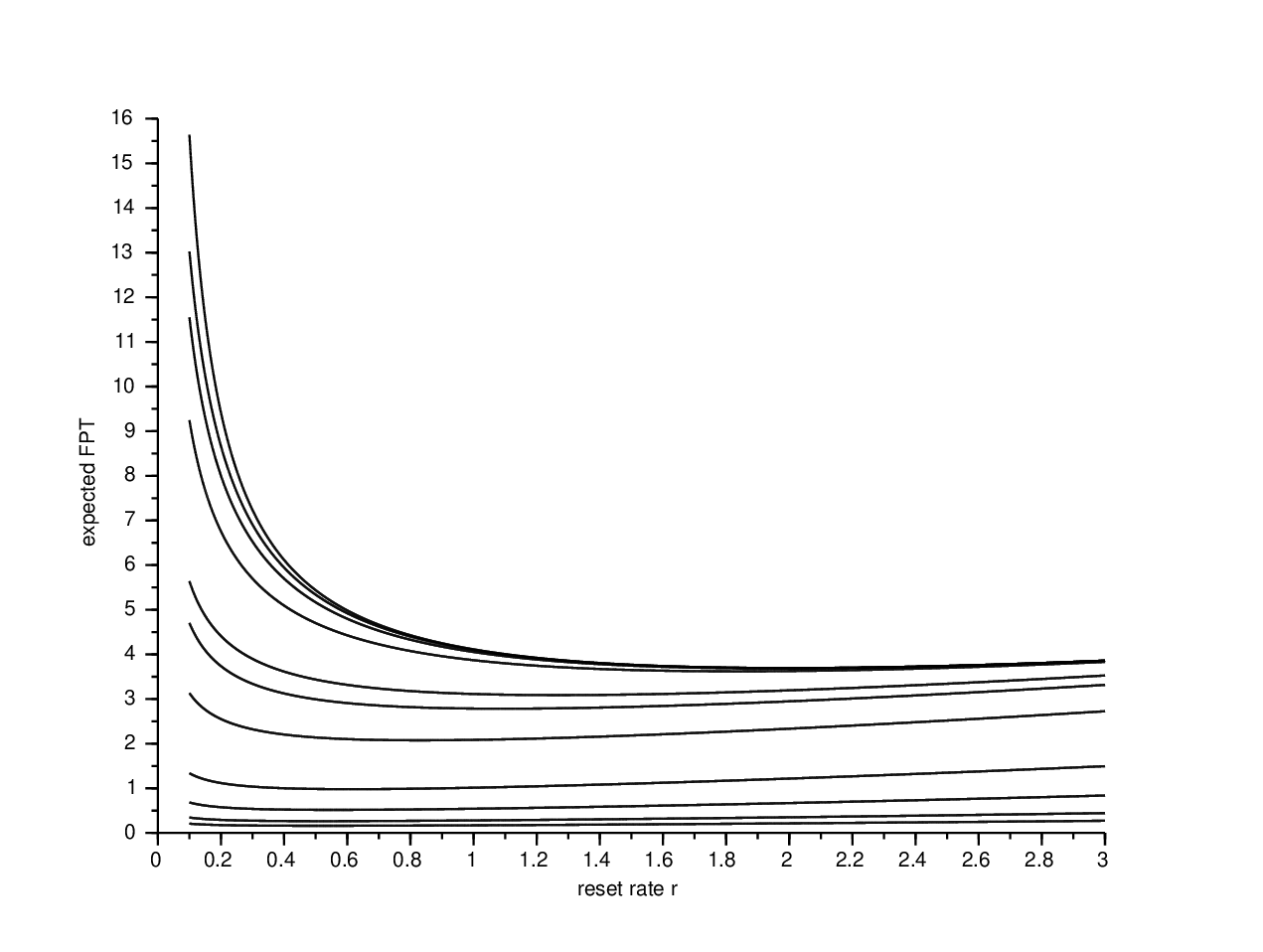}
\caption{For BM with resettings, the figure shows the graphs  of the expected FPT $T(x,r),$  as functions of $r>0,$ for fixed $x_R=1$ and for the values of $x$
contained in  the first column of Table \ref{tab2} (on the horizontal axes $r);$
the point of mimimum $r_{m}(x)$ increases from $\alpha = \frac 1 {2x_R^2} = \frac 1 2 $ (attained at  $x = 0.0001)$ to $\beta = \frac 2 {x_R^2} = 2$ (obtained for large $x >0).$
}
\label{FPTmin}
\end{figure}

Figure \ref{FPTmin} shows the graphs of $T(x,r) ,$ as functions of $r>0,$ for fixed $x_R=1,$ and for the values of $x$ contained in the first column of  Table \ref{tab2};
the greater the value of $x,$ the higher the corresponding curve, and the greater $r_{m}(x),$ namely, the abscissa of the point with horizontal tangent shifts more and more to the right.\par\noindent
It can be noted that, varying the curves,
the abscissa of the minimum, $r_{m}(x),$ increases from $ \alpha = \frac 1 {2 x_R^2}= \frac 1 2 ,$ obtained at $x = 0.0001 ,$ to $\beta = \frac 2 {x_R^2} = 2,$ obtained for large values of $x$
(one has $\alpha = \frac 1 2 $ and $\beta =2,$ being $x_R=1).$ \par
\begin{figure}[!h]
\centering
\includegraphics[height=0.35 \textheight]{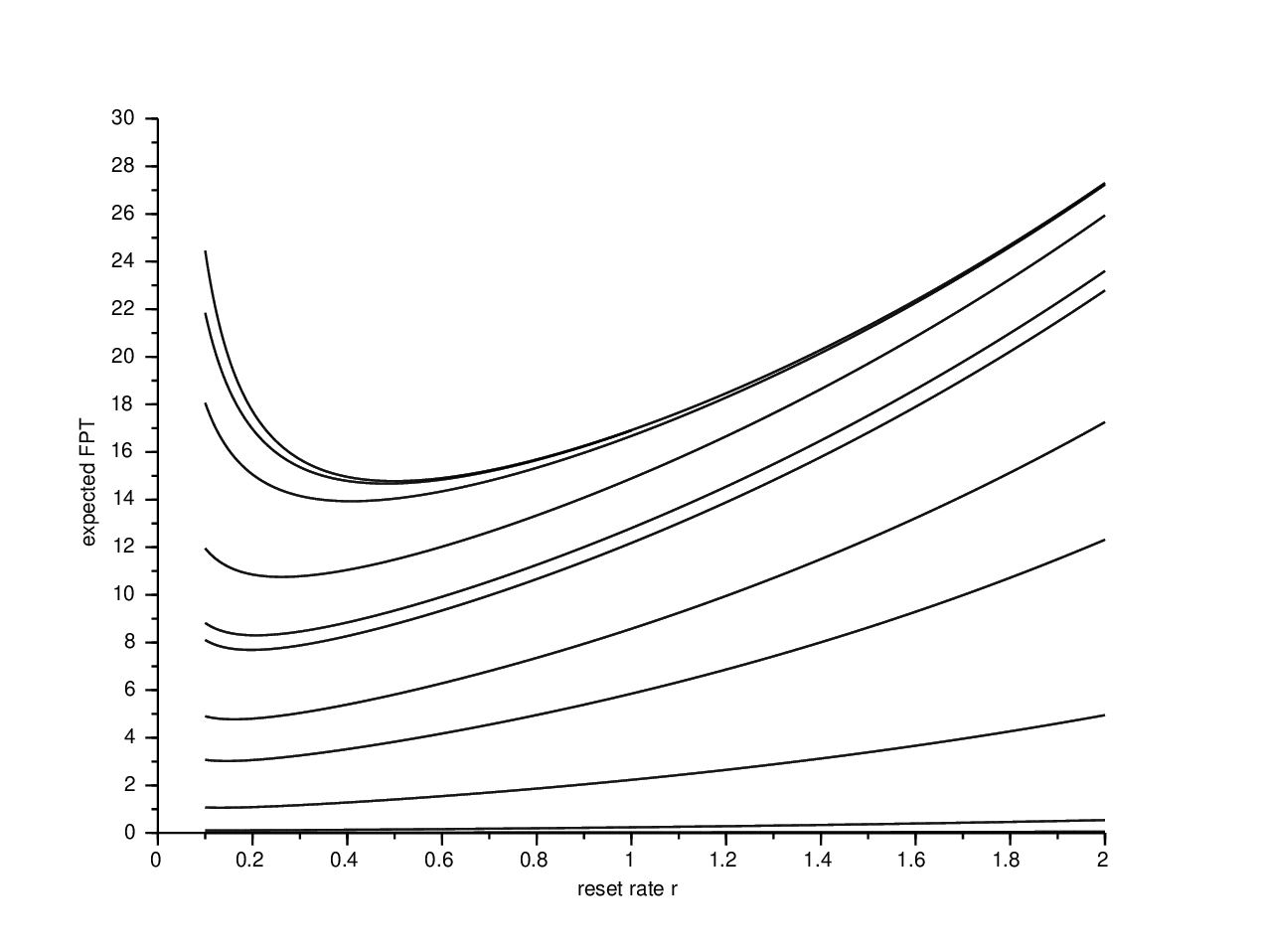}
\caption{For BM with resettings, the figure shows the graphs  of the expected FPT $T(x,r),$   as functions of $r>0,$ for fixed $x_R=2$ and for the values of $x$
contained in  the first column of Table \ref{tab2} (on the horizontal axes $r);$
the point of mimimum $r_{m}(x)$ increases from $\alpha = \frac 1 {2x_R^2} = 1/8$ (obtained at  $x =0.0001)$ to $\beta = \frac 2 {x_R^2} = 1/2$ (obtained for large $x >0).$
}
\label{FPTmin3}
\end{figure}

Actually, the qualitative behaviors of ${r_m(x)}$ and $m(x)$ do not appear to depend on the value of $x_R ; $ we have graphically shown them for $x_R=1,$ however they are similar, for other values of $x_R.$
In fact, in  Figure \ref{FPTmin3} the graphs of $T(x,r),$ as functions of $r,$  are reported for the same set of values of $x$  of  Figure \ref{FPTmin}, but for fixed $x_R=2:$ the qualitative behaviors of the curves are
similar to those in Figure \ref{FPTmin}. \par\noindent

\begin{figure}[!h]
\centering
\includegraphics[height=0.23 \textheight]{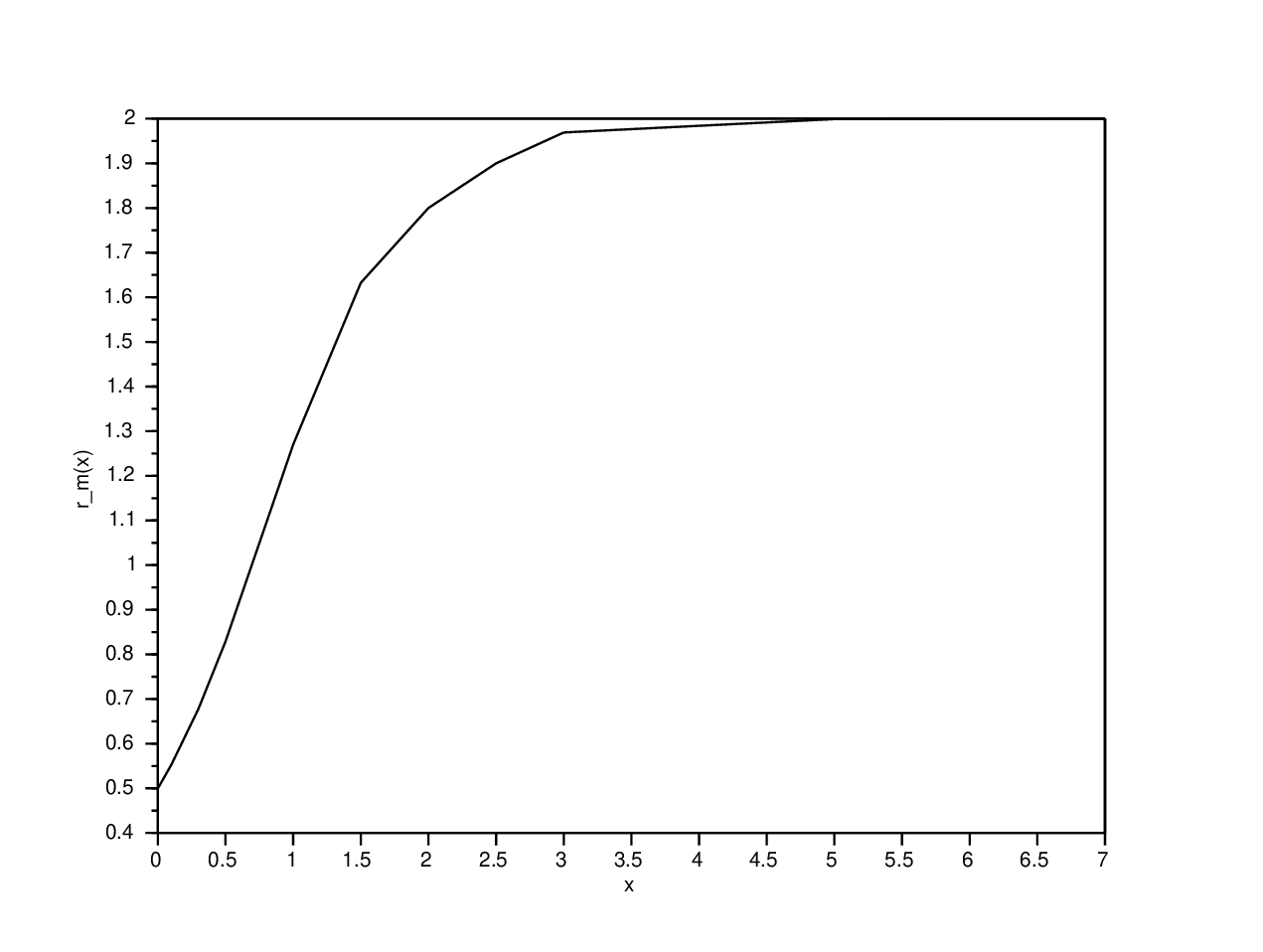}
\includegraphics[height=0.23 \textheight]{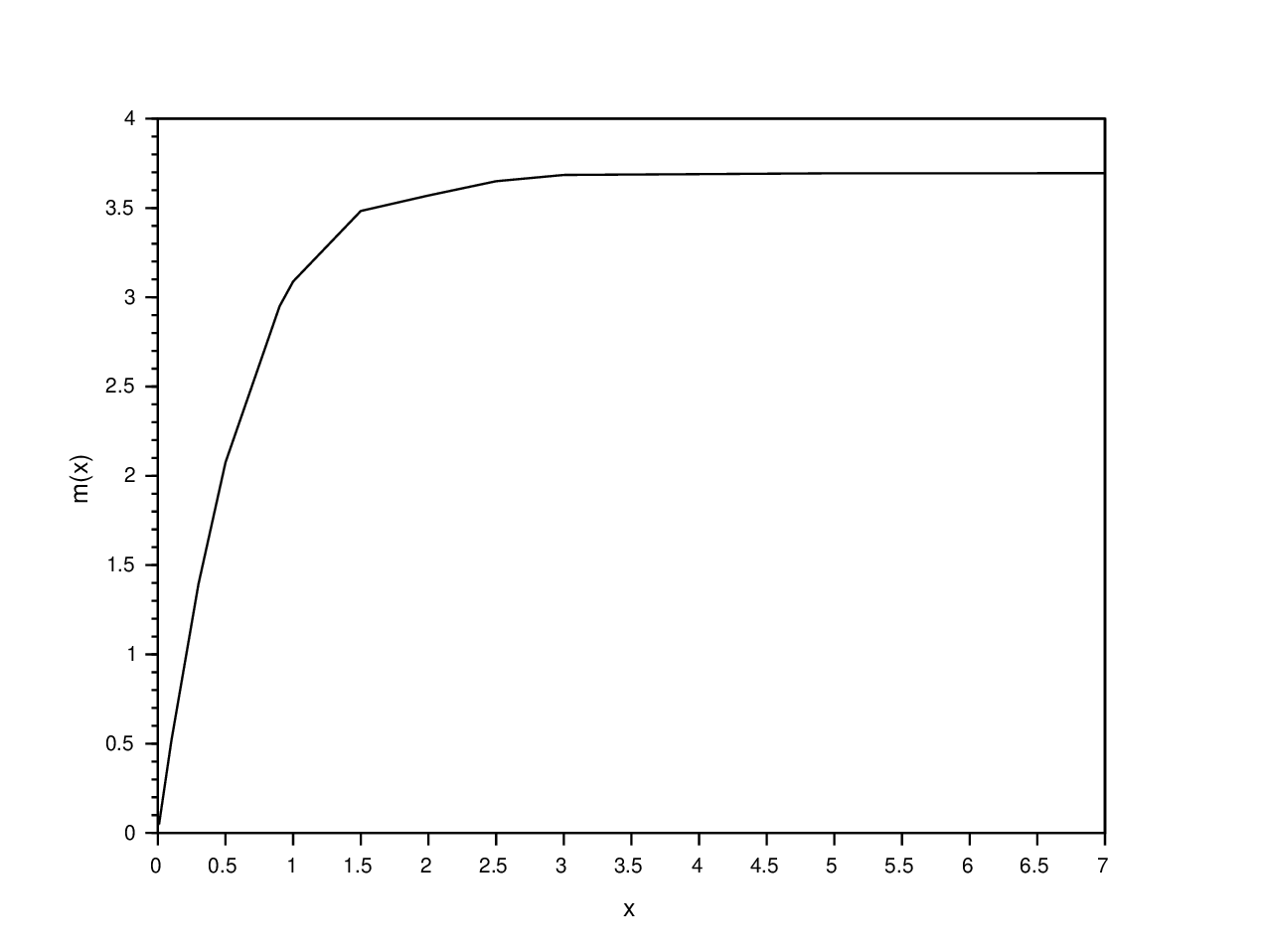}
\caption{For BM with resettings, the figure shows the graphs  of $r_{m}(x) $ (left panel), and
 $m(x)= T (x,r_{m}(x))$ (right panel), as functions of $x >0,$ for fixed $x_R=1$ and for the values of $x$ contained in  the first column of  Table \ref{tab2};
   here, $T(x,r)$ is the expected FPT.  The function  $r_{m}(x) $ increases from $\alpha = 1/2$ to $\beta =2,$  while
 $m (x)$ increases from about $0$ to  $\frac {e^ {x_R \sqrt {2\beta}}} \beta = \frac 1 2  e ^ 2 \approx 3.695 \ .$
}
\label{minEFPT}
\end{figure}

In Figure \ref{minEFPT} we report the
graphs of $r_{m}(x) $ (left panel), and
 $m(x)= \min _ {r \ge 0}T(x,r)$ (right panel), as  functions of $x >0,$ for fixed $x_R=1;$ note that $r_{m}(x) $ increases from $\alpha =  \frac 1 {2 x_R^2}= 1/2$ to
 $\beta = \frac 2 {x_R^2} =2$  (for large $x >0),$ while $m(x)$ increases from a value of about zero to  $\frac {e^ {x_R \sqrt {2\beta}}} \beta = \frac 1 2  e ^ 2 x_R ^2  = 3.695 \ .$

\bigskip

\noindent{\bf Remark}
If one keeps $r$ and $x$ fixed, then the expected FPT of ${\cal X}(t)$ through zero results to be an increasing function of $x_R >0,$ so its minimum is obtained for $x_R=0.$

\bigskip

\bigskip
\newpage

\noindent {\bf (b) Drifted BM with resetting} \par\noindent
For non-zero drift $\eta$ the expected FPT is given by \eqref{EFPTBMreset}.
\par\noindent
We proceed to minimize $T(x,r)= E[\tau (x,r)] $ with respect to the resetting rate $r,$ as done in the zero-drift case. As the results are in part similar to those concerning
the case $\eta =0,$ for brevity we limit ourselves to report them in tabular form, without graphical outputs. \par\noindent
Unlike the case of drifted BM without resetting $(r=0),$ the expected FPT, $T(x,r),$ results to be finite also for $\eta >0,$ for all $x >0$ and $r>0$ (see e.g. \cite{abundo:FPA2023}).
As before,
$
\lim _ {r \rightarrow 0^+} T(x,r) = + \infty, \ {\rm as \ well \ as } \ \lim _ {r \rightarrow + \infty } T( x,r) = + \infty.
$
For fixed reset position $x_R >0$ and starting point $x >0,$ $T(x,r),$  as a function of $r,$ attains the unique global minimum at a value
$
r_m(x)= arg \left ( \min _ {r \ge 0 } T (x,r) \right ).
$
\bigskip

\noindent {\bf (i) The case when $x=x_R $ } \par\noindent
Taking $x= x_R$ in \eqref{EFPTBMreset},  for $r>0$ one gets:
\begin{equation}
T(x,r)= \frac 1 r  \left ( e^{x \left (\eta + \sqrt {\eta ^2 +2r} \ \right )} -1 \right ).
\end{equation}
In Tables \ref{tabdriftBMmupm01} and \ref{tabdriftBMmupm1}, we report the values of  $r_m(x), \ m(x)$ and $T(x,0)$  as functions of $x,$ for  $\eta = \pm 0.1 $ and $\eta = \pm 1  ;$
as it must be, the numbers in the third column are not greater than those in the fourth column. \par\noindent
As we see, in all cases
$r_m(x)$ decreases, while $m(x)$ increases, as in the zero-drift case with $x=x_R$ (see Table \ref{tab1}). Of course, for $\eta = \pm 0.1,$ due to the small absolute value of the drift,  the values appear to be
closer to those regarding $\eta =0,$ than they do for $\eta = \pm 1 .$ \par\noindent
Note that, for $\eta = -1$, it results $r_m(x)=0$ for $x = x_R \ge 2.$ In fact, since $\eta$ is negative and sufficiently large in absolute value, the expected FPT is rather small  even for $r=0;$ indeed, resetting to a large $x_R,$
instead of accelerating the passage through zero, slows it down.


 \begin{table}[!h]
\centering
\begin{tabular}{cccc}
\hline
 $x$ & $r_{m}(x)$ & $m(x)$ & $T(x,0) = -\frac x \eta $ \\
\hline
0.1    & 185. & 0.030  & 1 \\
0.5    & 4.859 & 0.725 & 5 \\
1.     & 1.159 & 2.727 & 10\\
2.     & 0.261 & 9.659 & 20 \\
3.     & 0.103 & 19.291 & 30 \\
5.     & 0.027 & 42.547 & 50 \\
10.     & 0. & 100.  & 100\\
\hline
\end{tabular}
\ \ \ \ \ \ \
\begin{tabular}{cccc}
\hline
 $x$ & $r_{m}(x)$ & $m(x)$ & $T(x,0)$ \\
\hline
0.1    & 128.072. & 0.0312  & $+ \infty$ \\
0.5    & 5.296 & 0.822 & $+ \infty$ \\
1.     & 1.377 & 3.505 & $+ \infty$\\
2.     & 0.370 & 15.955 & $+ \infty$ \\
3.     & 0.176 & 40.949 & $+ \infty$ \\
5.     & 0.071 & 149.025 & $+ \infty$ \\
10.     & 0.022 & 1216.25  & $+ \infty$\\
\hline
\end{tabular}

\caption{Minimization of the expected FPT, $T(x,r),$ of drifted BM with resetting: the tables report the values of $r_{m}(x), \ m(x)=T(x,r_{m}(x))$ and $T(x,0),$ numerically obtained for some values of $x =x_R >0$ for $\eta = -0.1 $ (left panel) and $\eta =0.1$ (right panel) }
\label{tabdriftBMmupm01}
\end{table}
\bigskip

\begin{table}[!h]
\centering
\begin{tabular}{cccc}
\hline
 $x$ & $r_{m}(x)$ & $m(x)$ & $T(x,0) = -\frac x \eta $ \\
\hline
0.1    & 114.811 & 0.027  & 0.1 \\
0.5    & 2.744 & 0.425 & 0.5 \\
1.     & 1.0008 & 0.99 & 1\\
2.     & 0 & 2 & 2 \\
3.     & 0 & 3 & 3 \\
5.     & 0 & 5 & 5 \\
10.     & 0 & 10  & 10\\
\hline
\end{tabular}
\ \ \ \ \ \ \
\begin{tabular}{cccc}
\hline
 $x$ & $r_{m}(x)$ & $m(x)$ & $T(x,0)$ \\
\hline
0.1    & 137.767 & 0.0350  & $+ \infty$ \\
0.5    & 7.154 & 1.490 & $+ \infty$ \\
1.     & 2.272 & 12.162 & $+ \infty$\\
2.     & 0.802 & 231.053 & $+ \infty$ \\
3.     & 0.4620 & 2786.840 & $+ \infty$ \\
5.     & 0.2439 & 270967. & $+ \infty$ \\
10.     & 0.1104 & $125 \times 10^{10}$  & $+ \infty$\\
\hline
\end{tabular}

\caption{Minimization of the expected FPT, $T(x,r),$ of drifted BM with resetting: the tables report the  values of $r_{m}(x), \ m(x)=T(x,r_{m}(x))$ and $T(x,0),$ numerically obtained for some values of $x =x_R >0$ for $\eta = - 1 $ (left panel) and $\eta =1$ (right panel) }
\label{tabdriftBMmupm1}
\end{table}
\bigskip

\bigskip
\newpage

\noindent {\bf (ii) The case when $x \neq x_R$ } \par\noindent
Now,
$T(x,r)= \frac 1 r e^{x_R \left (\eta + \sqrt {\eta ^2 +2r} \ \right )} \left ( 1- e^{-x \left (\eta + \sqrt {\eta ^2 +2r} \ \right )} \right ).$ \par\noindent
In Tables \ref{tabdriftBMxnotxrmupm01} and \ref{tabdriftBMxnotxRmupm1}, we report the values of  $r_m(x), \ m(x)$ and $T(x,0)$  as functions of $x,$ for $x_R =1$ and some values of  $\eta .$ As we see,
$r_m(x)$ and $m(x)$ are both increasing function of $x,$ as in the zero-drift case with $x \neq x_R =1$ (see Table \ref{tab2}); as $x$ becomes large they appear to approximate certain numbers
$\alpha$ and $\beta, $ respectively, as in the zero-drift case.  Of course, for $\eta = \pm 0.1$ the values are
more similar to those regarding $\eta =0.$
\par
As in the case of zero drift, for other values of $x_R$ (not reported), we have obtained  similar behaviors.

 \begin{table}[!h]
\centering
\begin{tabular}{cccc}
\hline
 $x$ & $r_{m}(x)$ & $m(x)$ & $T(x,0) = -\frac x \eta $ \\
\hline
0.1    & 0.448 & 0.427  & 1 \\
0.5    & 0.708 & 1.777 & 5 \\
1.     & 1.1593 & 2.727 & 10\\
2.     & 1.7973 & 3.269 & 20 \\
3.     & 1.963 & 3.339 & 30 \\
5.     & 2.0035 & 3.351 & 50 \\
10.     & 2.0049 & 3.3513  & 100\\
50.     & 2.0049 & 3.3513  & 500\\
\hline
\end{tabular}
\ \ \ \ \ \ \
\begin{tabular}{cccc}
\hline
 $x$ & $r_{m}(x)$ & $m(x)$ & $T(x,0)$ \\
\hline
0.1    & 0.657 & 0.624  & $+ \infty$ \\
0.5    & 0.948 & 2.425 & $+ \infty$ \\
1.     & 1.377 & 3.505 & $+ \infty$\\
2.     & 1.873 & 4.028 & $+ \infty$ \\
3.     & 1.982 & 4.085 & $+ \infty$ \\
5.     & 2.0044 & 4.093 & $+ \infty$ \\
10.     & 2.0049 & 4.093  & $+ \infty$\\
50.     & 2.0049 & 4.093  & $+ \infty$\\
\hline
\end{tabular}

\caption{Minimization of the expected FPT, $T(x,r),$ of drifted BM with resetting: the tables report   the values of $r_{m}(x), \ m(x)=T(x,r_{m}(x))$ and $T(x,0),$ numerically obtained for some values of $x,$ for $x_R=1$ and $\eta = -0.1 $ (left panel) and $\eta =0.1$ (right panel). Note that, as $x $ becomes large, $r_m(x)$  and $m(x)$ approach the values $2.0049$ and $3.3513$ (respect.) for  $\eta = -0.1, $ while they approach the values $2.0049$ and $4.093$ (respect.) for  $\eta = 0.1 \ . $ }
\label{tabdriftBMxnotxrmupm01}
\end{table}
\bigskip

\begin{table}[!h]
\centering
\begin{tabular}{cccc}
\hline
 $x$ & $r_{m}(x)$ & $m(x)$ & $T(x,0) = -\frac x \eta $ \\
\hline
0.1    & 0 & 0.1  & 0.1 \\
0.5    & 0  & 0.5 & 0.5 \\
1.     & 0.0011 & 1. & 1\\
2.     & 1.607 & 1.566 & 2 \\
3.     & 2.193 & 1.675 & 3 \\
5.     & 2.396 & 1.702 & 5 \\
10.     & 2.4141 & 1.703  & 10\\
50.     & 2.4142 & 1.703  & 10\\
\hline
\end{tabular}
\ \ \ \ \ \ \
\begin{tabular}{cccc}
\hline
 $x$ & $r_{m}(x)$ & $m(x)$ & $T(x,0)$ \\
\hline
0.1    & 1.6009 & 3.466  & $+ \infty$ \\
0.5    & 1.9814 & 10.194 & $+ \infty$ \\
1.     & 2.272 & 12.162 & $+ \infty$\\
2.     & 2.405 & 12.575 & $+ \infty$ \\
3.     & 2.413 & 12.588 & $+ \infty$ \\
5.     & 2.4142 & 12.589 & $+ \infty$ \\
10.     & 2.4142 & 12.589 & $+ \infty$\\
50.     & 2.4142 & 12.589 & $+ \infty$\\
\hline
\end{tabular}

\caption{Minimization of the expected FPT, $T(x,r),$ of drifted BM with resetting: the tables report the  values of $r_{m}(x), \ m(x)=T(x,r_{m}(x))$ and $T(x,0),$ numerically obtained for some values of $x,$ for $x_R=1$ and $\eta = -1 $ (left panel) and $\eta = 1$ (right panel). Note that, as $x$ becomes large, $r_m(x)$  and $m(x)$ approach the values $2.4142$ and $1.703$ (respect.) for  $\eta = -1, $ while they approach the values $2.4142$ and $12.589$ (respect.) for  $\eta = 1 \ . $   }
\label{tabdriftBMxnotxRmupm1}
\end{table}
\bigskip

\subsubsection{Optimization of the FET}
\noindent {\bf (a) BM with resetting ($\eta =0)$} \par\noindent
Let $\mathcal X (t)$ be  BM with resetting. For $x_R $ and  $x \in (0,b), \ \tau  (x,r)$ represents now
the first-exit time (FET) of $\mathcal X  (t)$ from the interval $(0,b),$ under the condition that $\mathcal X(0)=x ,$ namely:
\begin{equation}
\tau  (x,r)= \min \{t>0: \mathcal X(t) \notin (0,b) | \mathcal X(0) =x  \}, \ x \in (0,b).
\end{equation}
\par\noindent
By taking $\eta =0$ in formula  \eqref{EtauBMdriftreset}, we obtain the expected FET:
\begin{equation} \label{EFET}
 T(x,r):= E[\tau  (x,r)] = \frac {\sinh \left (b \sqrt {2r} \ \right )- \sinh \left (x \sqrt {2r} \ \right ) - \sinh \left ((b -x) \sqrt {2r} \ \right )   }  {r \left [\sinh \left (x_R \sqrt {2r} \ \right ) + \sinh \left ((b - x_R) \sqrt {2r}  \ \right ) \right ]}.
\end{equation}
(the notation includes the dependence on $x$ and $r,$ but not on $x_R,$ for simplicity). \par\noindent
For  $r \ge 0,$ one has $T(0,r)= T(b,r) =0. $ \par\noindent
Since, for fixed $r, \ T(x,r)$ is symmetric with respect to $x=b/2,$ that is,
$
T(x,r)= T(b-x, r), \ x \in (0,b) ,
$
one has to calculate $T(x,r)$ only for $x \in (0, b/2].$ \par\noindent
Even though we consider the boundaries $0$ and $b,$  an analogous formula for $T(x,r)$ can be obtained for any couple of boundaries $a, \ b$  (see \cite{abundo:MCAP2025}). \par\noindent
Unlike the previous case of one boundary (namely the FPT), the expected FET of $\mathcal X  (t)$ from $(0,b)$ is finite also for $r=0,$ that is, for BM without resetting; in fact,  in this case one has $T(x,0)= x(b-x).$ \par
As easily seen,  for fixed reset position $x_R \in (0,b)$ and starting point $x \in (0,b),$ the expected FET, $T( x,r) ,$ as a function of $r,$ attains the unique global minimum at the value
$
r_m(x)= arg \left ( \min _ {r \ge 0 } T (x,r) \right ) .
$
As in Section 5.1.1, for fixed $x_R \in (0,b),$ we will  find  the optimal reset value $r_{m}(x)$ and the value of
the minimum, $m(x)= T ( x,r_{m}(x)) ,$ as functions of $x \in (0,b).$ \par\noindent
Note that, if one keeps $r$ and $x$ fixed, then the expected FET, as a function of $x_R \in (0,b),$ has its global maximum at $x_R = b/2,$ while the minimum is attained
at $x_R=0$ or $x_R=b.$ \bigskip

\noindent {\bf (i) The case when $x=x_R$ } \par\noindent
For $x = x_R \in (0,b), \ r >0 ,$ Eq. \eqref{EFET} becomes:
\begin{equation}
T(x,r)= \frac {\sinh (b \sqrt {2r} )} {r [\sinh((b-x) \sqrt {2r}) + \sinh(x \sqrt {2r})]} - \frac 1 r .
\end{equation}
For small  $x \in (0,b)$ the above quantity, as a function of $r >0,$ is first decreasing and then increasing, hence there exists a global minimum point $r_m(x) >0.$ For larger $x$ instead (specifically for $x \ge 0.2 ,$ when $b=1), \ T(x,r)$ is increasing, so $r_m(x)=0.$
\par\noindent
In Table \ref{tab2b}, we report the values of $r_{m}(x), \ m(x) = \min _ {r \ge 0 } T (x,r)= T(x, r_{m}(x)),$ and $T(x,0)= x(1-x),$ numerically obtained  for $b=1,$ for some values of
$x = x_R$ belonging to  $(0,1/2].$ For $x \in (1/2, 1)$ they are obtained by symmetry, being $r_m(x)= r_m(1-x)$ and $m(x) = m(1-x).$
As before, the numbers in the third column are  not greater than those in the fourth column, which provide the expected FPT, i.e. $T(x,0)= x(1-x),$ in the no-resetting case
\par\noindent

\begin{table}[!h]
\begin{center}
\begin{tabular}{cccc}
\hline
 $x$ & $r_{m}(x)$ & $m(x)$ & T(x,0) \\
\hline
0.1    & 126.972 & 0.0308 & 0.09 \\
0.2    & 28.442 & 0.1221 & 0.16 \\
0.25    & 10.131 & 0.1795 & 0.1875 \\
0.27    & 2.610 & 0.1965 & 0.1971 \\
0.275   & 0.580 & 0.199 & 0.1993 \\
0.28    & 0 & 0.201 & 0.2016 \\
0.3    & 0 & 0.21 & 0.21 \\
0.4    & 0 & 0.24  & 0.24\\
0.5    & 0 & 0.25  & 0.25\\
\hline
\end{tabular}
\end{center}
\caption{For BM with resetting the table reports the values of $r_{m}(x), \ m(x),$ and $T(x,0),$  numerically obtained in the two-boundary case with $b=1,$ for some values of
$x = x_R$ belonging to $(0,1/2]$ (for $x \in (1/2, 1)$ they are obtained by symmetry, being $r_m(x)= r_m(1-x)$ and $m(x) = m(1-x).$  Here, $T(x,r)$ is the expected FET. }
\label{tab2b}
\end{table}

\bigskip

\noindent {\bf (ii) The case $x \neq x_R $ } \par\noindent
In this case, we observe a more intricate scenario.
As an example, we report in Table \ref{tab4} the values of $r_{m}(x)$ and $m(x) = \min _ {r \ge 0 } T (x,r)= T(x, r_{m}(x)),$ numerically obtained  by
the secant method, for some values of $x \in (0,b/2],$ with $b=1$ and $x_R =0.2$ (the values for $x \in (b/2, 1)$ can be obtained by symmetry). \par\noindent
Actually, $r_m(x)$ increases for $x \in (0, 1/2),$ whereas it decreases for $x \in (1/2, 1),$ namely $r_m(x)$ attains its maximum at $x = 1/2;$
the same thing happens for $m(x).$ \par\noindent
Even now,  all the values of $m(x)$  in Table \ref{tab4} are smaller than $T(x,0)= x(1-x),$ which are the expected FETs in the no-resetting case, thus the choice of the resetting rate $r= r_m(x)$
given in the second column expedites the FET.
 \par
\begin{figure}[!h]
\centering
\includegraphics[height=0.35 \textheight]{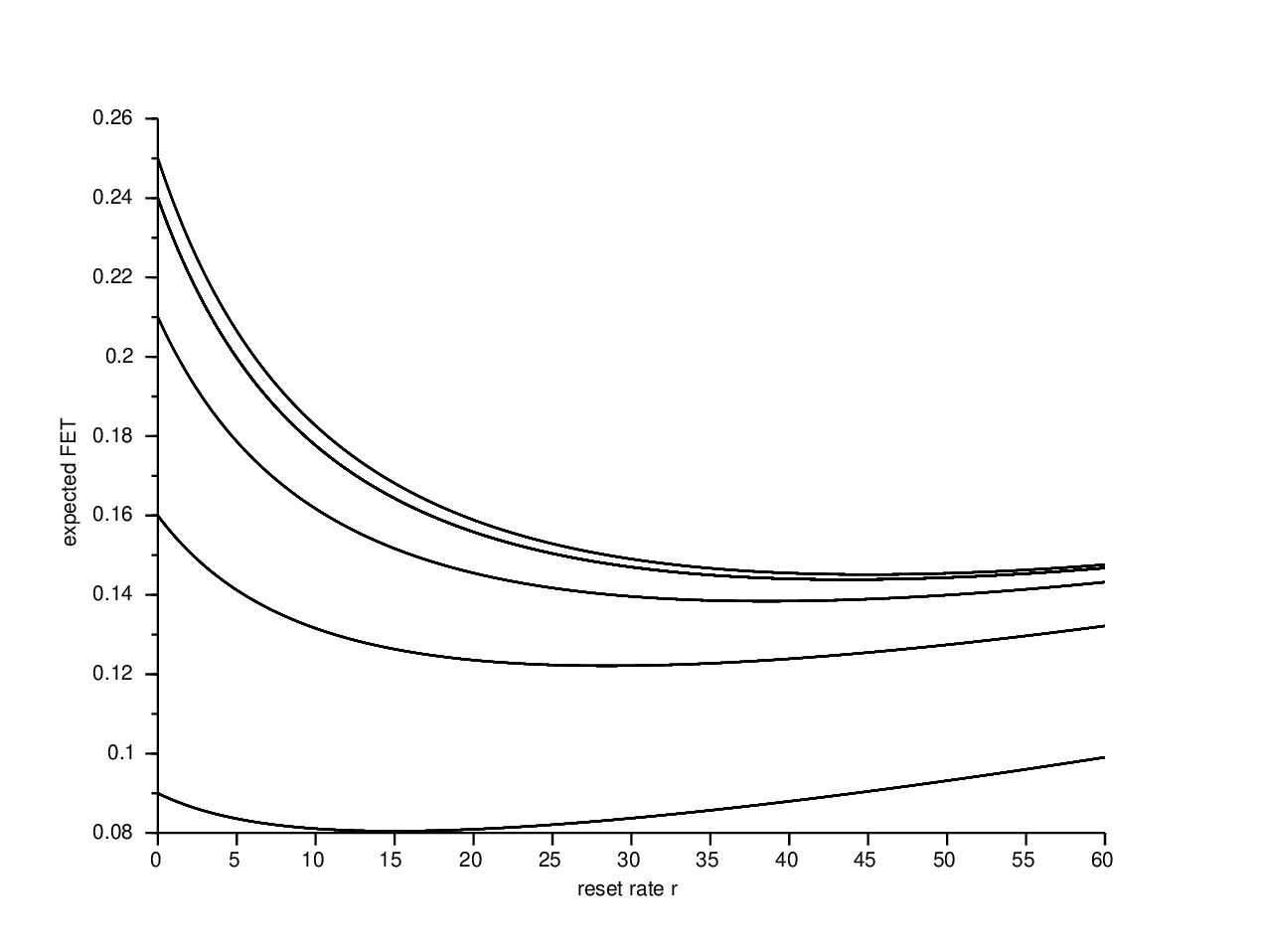}
\caption{For BM with resetting, the figure shows the graphs of the expected FET $T(x,r),$  as functions of $r>0,$ for $b=1, \ x_R=0.2 \ , $  and the values of $x$ going from $0.1$ to $0.5,$ with step $0.1$
 (on the horizontal axes $r);$  the lower and upper curve correspond to $x = 0.1$ and  $x = 0.5,$ respectively. As $x$ increases from $0.1$ to $0.5,$ the value
 $r_{m}(x)$ at which the minimum of $T(x,r)$ is attained, increases from $14.948$ to $45.009 ,$ while $m(x)= T(x, r_m(x))$ increases from $0.0804$ to $0.1451 $ (see Table \ref{tab4}).
}
\label{FETminDouble}
\end{figure}

In Figure \ref{FETminDouble}  the graphs of $T(x,r) ,$ as functions of $r>0,$  are reported  for the values of $x$ going from $0.1$ to $0.5,$ with step $0.1,$ for fixed $x_R=0.2 \ .$
The graphs of $T(x,r)$ for  $1/2 < x <1 ,$ can be obtained by symmetry, because
$T(x,r)= T(1-x, r).$
 The lower and upper curve correspond to $x = 0.1$ and  $x = 0.5,$ respectively; for  increasing values of $x ,$ the corresponding curves become higher and higher, and the abscissa of the minimum shifts more and more to the right. \par
By  evaluating numerically  $\alpha  := \lim _ {x \rightarrow 0 ^+} r_m(x)$ and $\beta :=  \max _ {x \in [0,1]} r_m(x),$ for $b=1$ and $x_R = 0.2 ,$ we  obtained the values $\alpha = 3.4325 ,$
and  $\beta =  45.009 ,$  while
$\gamma  :=  \lim _ {x \rightarrow 0 ^+} m(x)$ resulted to be approximately zero, and $\max _ {x \in [0,1]} m(x) = 0.1451$ (see Table \ref{tab4}). \par
\begin{figure}[!h]
\centering
\includegraphics[height=0.23 \textheight]{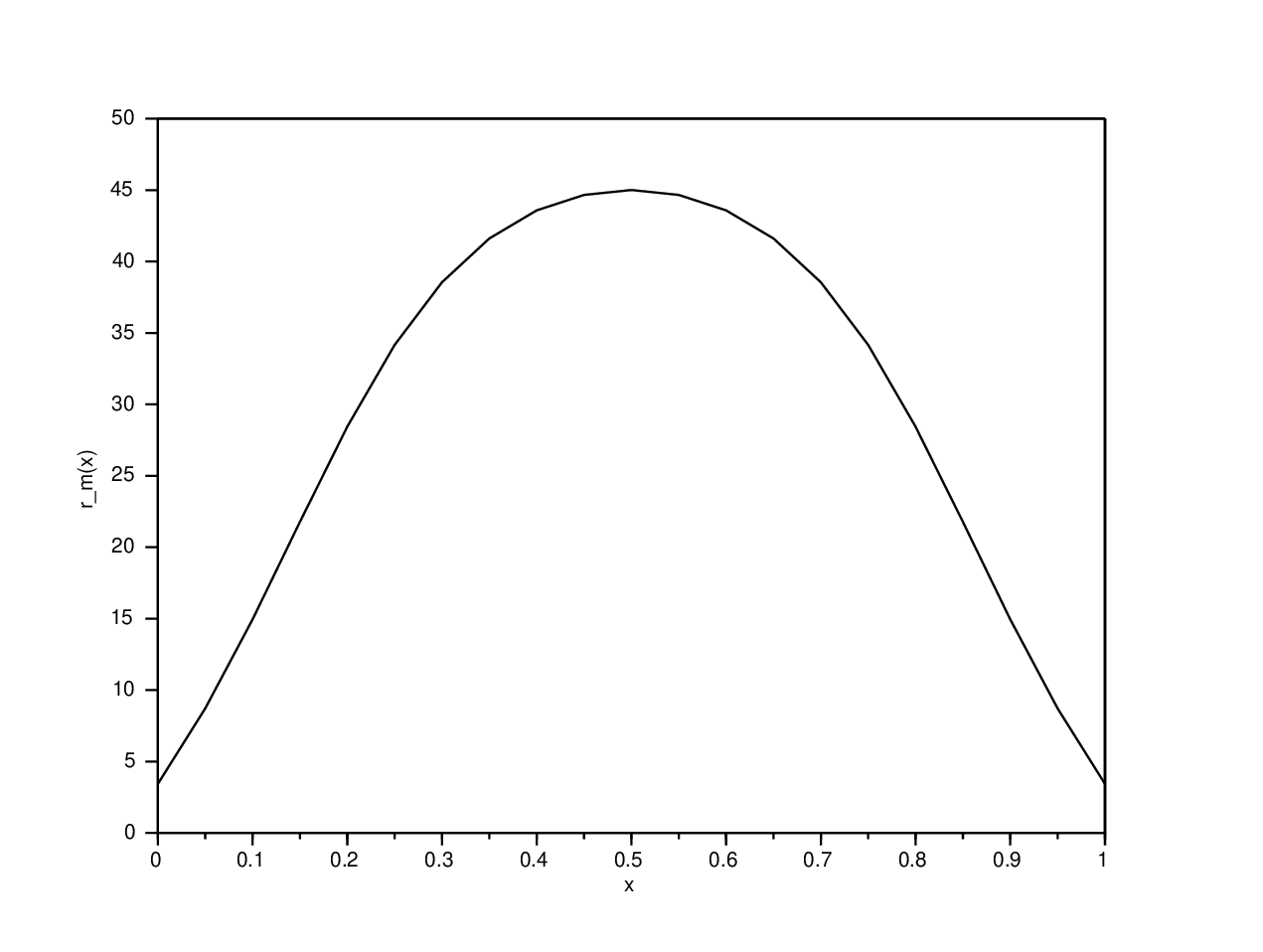}
\includegraphics[height=0.23 \textheight]{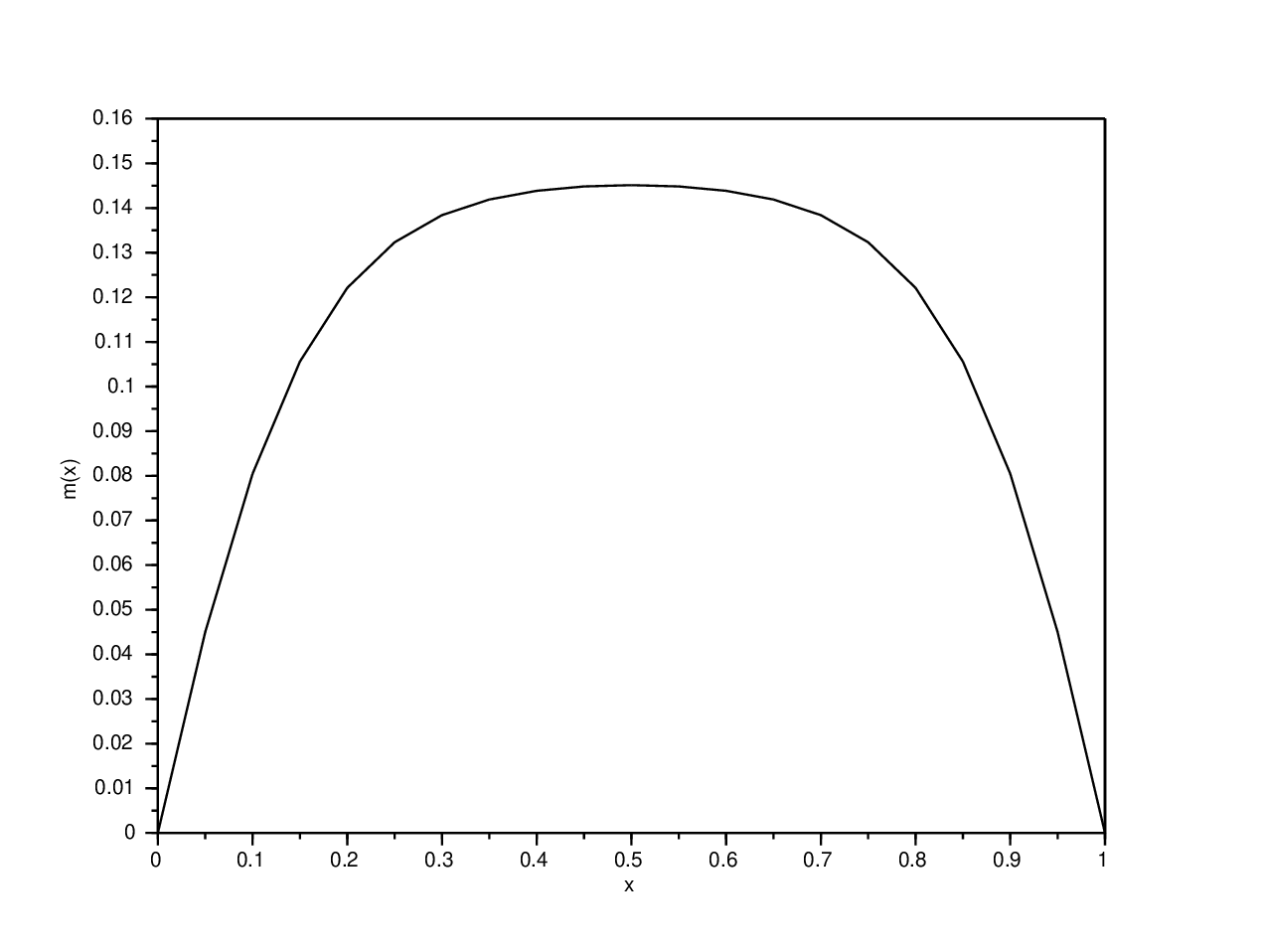}
\caption{For BM with resetting, the figure shows the graphs of $r_{m}(x) $ (left panel), and
 $m(x)= T (x,r_{m}(x))$ (right panel), as a function of $x \in (0,1),$ for   $x_R=0.2 \ ;$  here, $T(x,r)$ is the expected FET. The function  $r_{m}(x) $ increases from $3.4325$ (at $x=0$) to $45.009$ (at $x = 1/2),$ and then decreases;
similarly, $m(x)$ increases from $0$ to $0.1451 \ .$
}
\label{minEFPTDouble}
\end{figure}

In  Figure \ref{minEFPTDouble} we report the
graphs of $r_{m}(x) $ (left panel), and
 $m(x)= \min _ {r \ge 0}T(x,r)$ (right panel), as a function of $x \in (0,1),$ for $x_R=0.2;$ the values of $x$
 go from $x \simeq 0$ to $x \simeq 1$ (the values are taken from  Table \ref{tab4}). We see that $r_{m}(x) $ increases from $\alpha = 3.4325$ (at  $x \simeq 0$) to $\beta = 45.009$ (at $x = 1/2),$ after that it decreases, approximating  once again $\alpha$ at $x \simeq 1,$ while $m(x)$ increases from about zero to  $0.1451 ,$ after that it decreases again to about zero.
\par\noindent
As for the FPT, we have observed that $\alpha$ decreases as $x_R$ increases, while $\gamma = \lim _ {x \rightarrow 0} m(x) $ remains always small. \par
Actually, our computations showed that the qualitative behaviors of ${r_m(x)}$ and $m(x)$ do not depend on the value of $b$ for fixed  $x_R \in (0,b):$  they are similar, for any values of $b.$ \par
As a further example, we have reported in  Table \ref{tab5} the values of $r_{m}(x)$ and $m(x)$ numerically obtained for $b=1, \ x_R=0.3 ,$ and for the values of $x$ from $0$ to $1$ with step
$0.1 \ .$
We see that $r_m(x)$ remains always zero for these $x;$ in fact, the values of $m(x)$ in the third column are exactly the same ones as $x(1-x),$ which correspond to the FET in the no-resetting case $(r=0).$ For larger values of
$x_R $ we have observed the same thing. \par\noindent
However, unlike the case $x_R \le 0.2 ,$ there exists a  special value $\bar x_R \in (0.2, 0.3)$  at which we detected
a different behavior; this can be explained, by considering that, as $x_R$ becomes close enough to $1/2 ,$ also small non-zero values of the resetting rate $r$ imply large values of
the expected FET. In fact, if the process is reset  to a position close to $1/2,$ it takes more time to reach one of the ends of the interval $(0, 1),$ hence the minimum of the expected FET is obtained at $r=0.$
Of course, this happens for any values of $b,$ when $x_R$ is close to $b/2.$ \par\noindent
Since $r_m(x)$ is an increasing function of $x$ for fixed $x_R \in (0, 1/2),$  the feature of $r_m(x)$ can be captured  by computing
its maximum value (attained at $x=0.5)$ and  its minimum value (attained at $x \simeq 0).$
We  obtained that the
value of $x_R$ at which
the maximum value of $r_m(x)$ becomes  zero  (remaining zero for all $x)$  is approximately $0.295.$\par\noindent
We concluded that at $\bar x_R = 0.295$ there is a
transition, namely a
change of behavior with respect to the case when $x_R  \le 0.2;$ in fact, for $x > \bar x_R$ even the maximum value of $r_m(x)$ becomes zero.
\par
\begin{figure}[!h]
\centering
\includegraphics[height=0.23 \textheight]{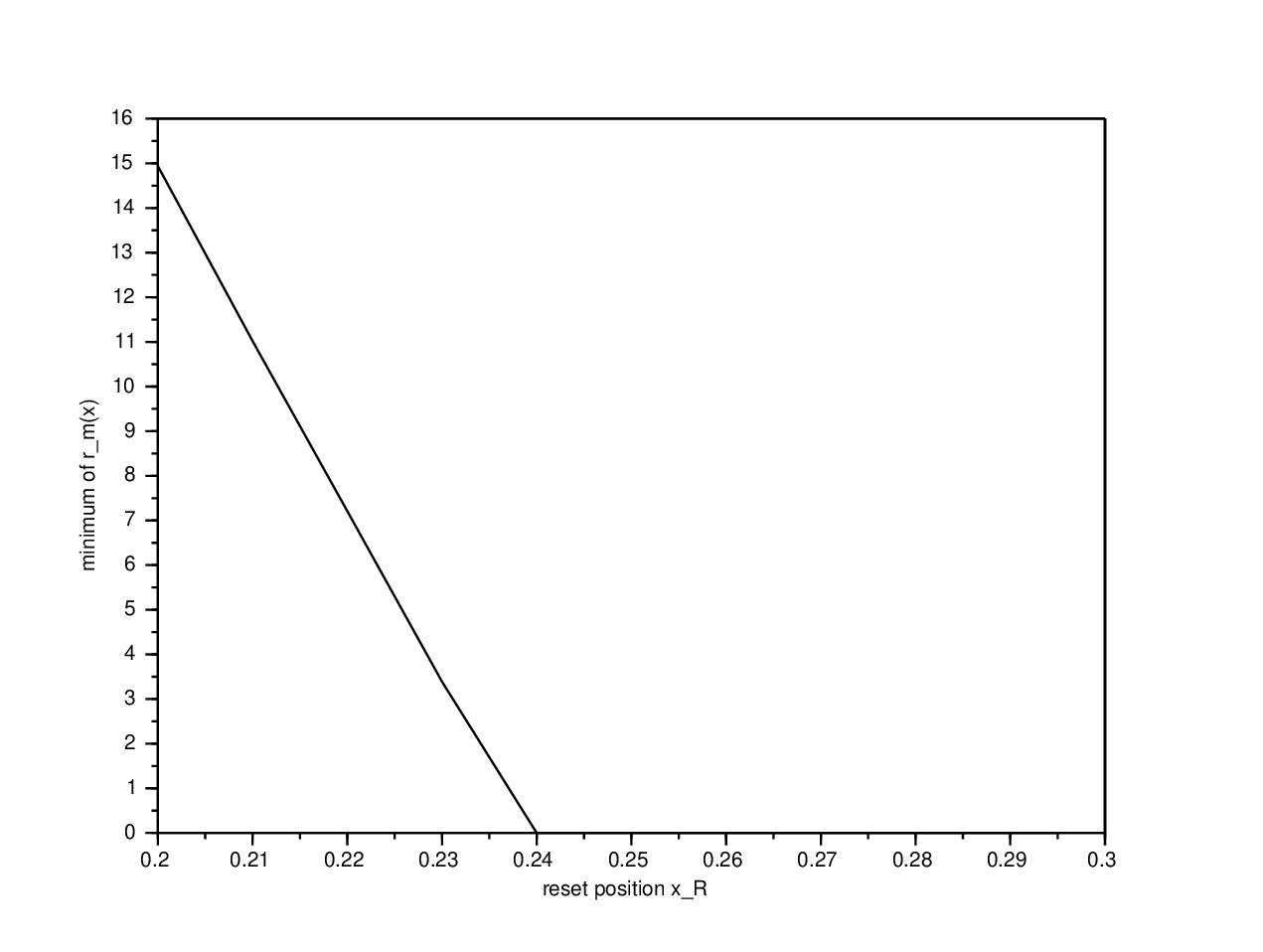}
\includegraphics[height=0.23 \textheight]{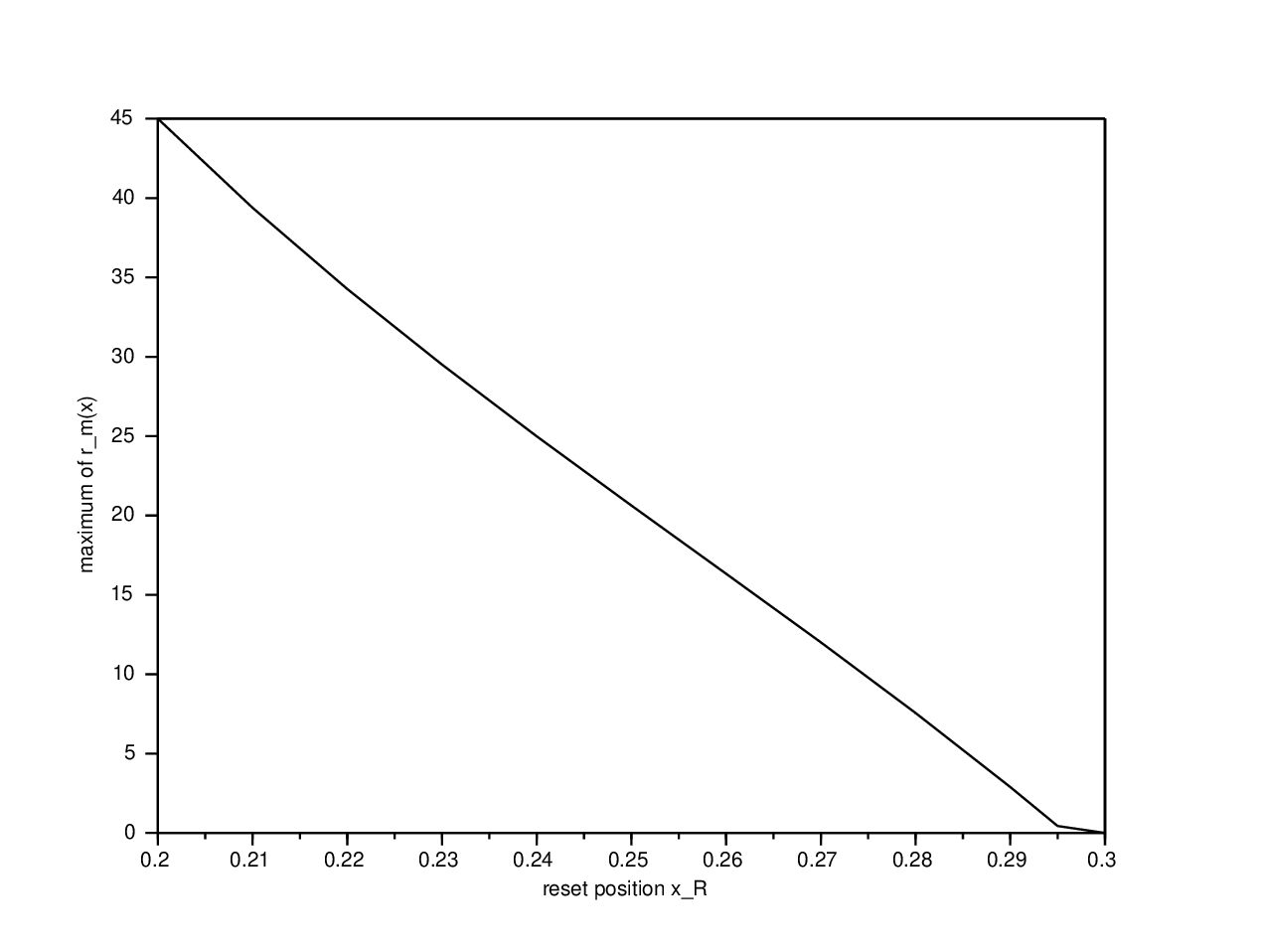}
\caption{For BM with resetting the figures show the graph of $r_m(0.1),$ i.e. the mimimum value of $r_m(x)$ for $x \in (0,1),$ as a function of $x_R \in [0.2, 0.3]$ in the two-boundary case for $b=1$ (left panel);
the graph of $r_m(0.5),$ namely the maximum value of $r_m(x),$ as a function of $x_R \in [0.2, 0.3]$ (right panel).
}
\label{EFPTDoubleTER.eps}
\end{figure}

In Figure \ref{EFPTDoubleTER.eps}  we show  the graph of $r_m(0.1),$ i.e. the minimum value of $r_m(x)$ for  $x \in (0.1,1),$  as a function of $x_R \in [0.2, 0.3]$ (left panel); we see that it decreases approximately linearly from $14.94$ to zero. In the right panel, we
show the graph of $r_m(0.5),$ namely the maximum value of $r_m(x),$ as a function of $x_R \in [0.2, 0.3];$ also it decreases approximately linearly  from $45.00$ to zero.
\par
\begin{figure}[!h]
\centering
\includegraphics[height=0.35 \textheight]{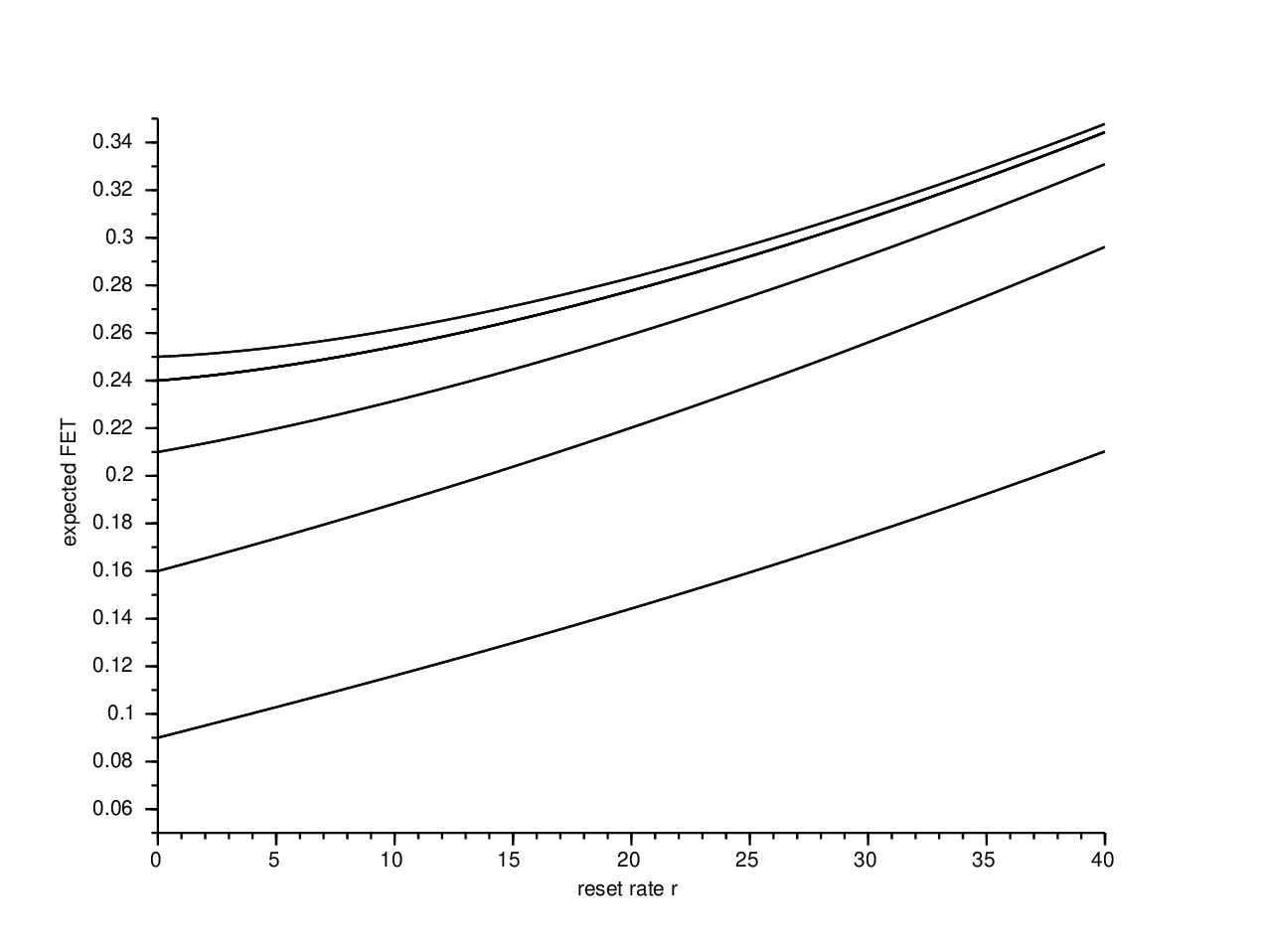}
\caption{For BM with resetting the figure shows  the graphs of the expected FET, $T(x,r),$  as functions of $r>0,$
for $b=1$ and $x_R=0.3 \ , $  for the values of $x$ going from $0.1$ to $0.5,$ with step $0.1$
 (on the horizontal axes $r);$  the lower and upper curve correspond to $x = 0.1$ and  $x = 0.5,$ respectively. As $x$ increases from $0.1$ to $0.5,$
 $r_{m}(x)$  remains always zero, while $m(x)= T(x, r_m(x))$ increases from $0.09$ to $0.25 .$
}
\label{FETminDoublexr03}
\end{figure}

In Figure \ref{FETminDoublexr03} we report the graphs of the expected FET $T(x,r) ,$  as functions of $r>0,$ for  $b=1, \ x_R=0.3 \ , $  and the values of $x$ going from $0.1$ to $0.5,$ with step $0.1 ;$
the lower and upper curve correspond to $x = 0.1$ and  $x = 0.5,$ respectively. As $x$ increases from $0.1$ to $0.5,$
 $r_{m}(x)$  remains always zero, while $m(x)=T(x, r_m(x))$ increases from $0.09$ to $0.25 $ (see Table \ref{tab5}). \par
\begin{figure}[!h]
\centering
\includegraphics[height=0.23 \textheight]{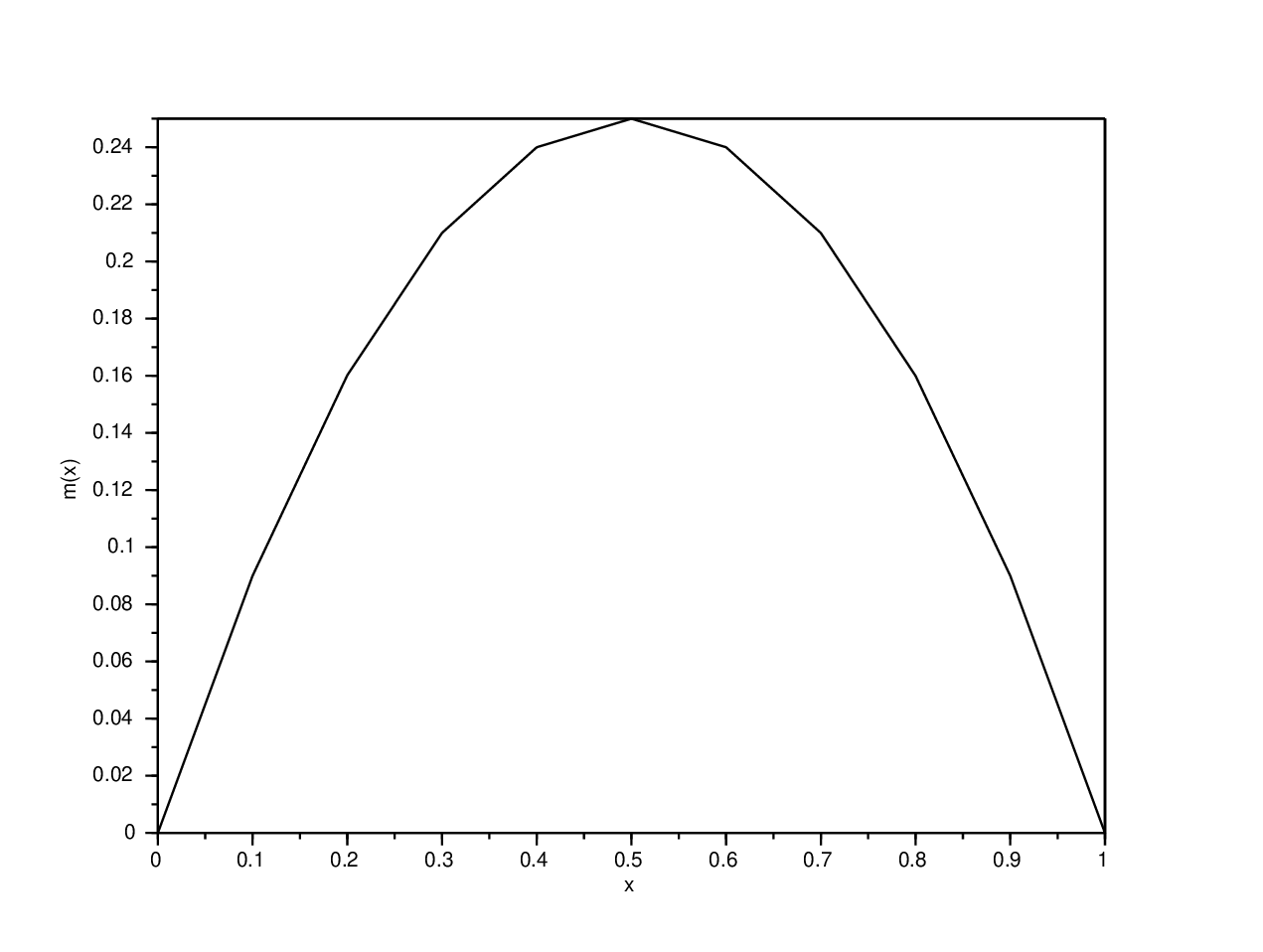}
\caption{For BM with resetting  it is shown the graph of $m(x)= \min _ {r \ge 0}T(x,r) $ in the two-boundary case for $b=1$ and $x_R=0.3 \ ,$ as a function of  $x \in (0, 1)$
 (on the horizontal axes $x).$ Here, $T(x,r)$ is the expected FET.   As $x$ increases
 $r_{m}(x)$  is always zero, while $m(x)$ increases from $0.09$ to $0.25 ,$ for $x \in [0, 1/2],$ coming back to the value
$0.09$ at $x=1.$
}
\label{EFPTDoubleBIS.eps}
\end{figure}

In  Figure \ref{EFPTDoubleBIS.eps} we report the graph of $m(x)= T(x,r_m(x)) ,$  as a function of $x \in (0,1),$ for $b=1$ and $x_R=0.3  \ .$
As $x$ increases,
$m(x)$ increases from $0.09$ to $0.25 $ (see  Table \ref{tab5}), while  $r_{m}(x)$  remains always zero (graph not reported).
For $x_R \in (0.3, 0.5),$ a similar situation can be  observed.\par
If  $x_R \in (0.2, \bar x_R),$ with  $ \bar x_R = 0.295,$ then the numerical computations  showed that   the functions $r_m(x)$ and $m(x)$  were not exactly zero, but they were very small; in particular,  $\alpha =
\lim _{x \rightarrow 0^+} r_m(x) $ and $ T(0, \alpha)$ appeared to be decreasing functions of the reset position $x_R,$ as in the one-boundary case. \bigskip

Actually, our computations confirmed all the particulars of the above scenario, for any value of $b,$ with a certain transition value $\bar x_R,$ close enough to $b/2.$

\begin{table}[!h]
\begin{center}
\begin{tabular}{cccc}
\hline
 $x$ & $r_{m}(x)$ & $m(x)$ & $T(x,0)$  \\
\hline
 $10^{-6} $ & 3.4325 & $9.8 \times 10^ {-8}$  & $ \simeq 0 $ \\
0.1    & 14.948 & 0.0804 & 0.09 \\
0.2    & 28.444 & 0.1221 & 0.16 \\
0.3    & 38.548 & 0.1384 & 0.21 \\
0.4    & 43.583 & 0.1438 & 0.24 \\
0.5    & 45.009 & 0.1451 & 0.25 \\
\hline
\end{tabular}
\end{center}
\caption{For BM with resetting, the table reports the values of $r_{m}(x), \ m(x)$ and $T(x,0),$ numerically obtained for  $b=1$ and $x_R=0.2,$ for some values of
$x \in (0,1/2]$ (the values for $x \in (1/2,1)$ can be obtained by symmetry with respect to $x= 1/2).$ Here, $T(x,r)$ is the expected FET.}
\label{tab4}
\end{table}

\par\noindent

\begin{table}[!h]
\begin{center}
\begin{tabular}{cccc}
\hline
 $x$ & $r_{m}(x)$ & $m(x)$ & $T(X,0)$ \\
\hline
0    & 0 & 0 & 0 \\
0.1    & 0. & 0.09 & 0.09 \\
0.2    & 0. & 0.16 & 0.16 \\
0.3    & 0. & 0.21 & 0.21 \\
0.4    & 0. & 0.24 & 0.24 \\
0.5    & 0. & 0.25 & 0.25 \\
\hline
\end{tabular}
\end{center}
\caption{For BM with resetting the table reports the values of $r_{m}(x), \ m(x)$ and $T(x,0),$ numerically obtained for $b=1$ and $x_R=0.3,$ for some values of $x \in [0,1/2];$ the values for $x \in (1/2,1)$ can be obtained by symmetry. Here, $T(x,r)$ is the expected FET. }
\label{tab5}
\end{table}

\noindent {\bf (b) Drifted BM with resetting ($\eta \neq 0$ )} \par\noindent
For non-zero drift $\eta$ the expected FET, $T(x,r),$  is (see \eqref{EtauBMdriftreset}):
$$ E[\tau (x,r)] = $$
\begin{equation}
= \frac {e^ {-b \eta }\sinh \left (b \sqrt {\eta ^2 +2r} \ \right )- e^ {- \eta x}\sinh \left (x \sqrt {\eta ^2 +2r} \ \right )
- e^ {-(b+x) \eta}\sinh \left ((b -x) \sqrt {\eta ^2 +2r} \ \right )   }  {r \left [ e^ {-x_R \eta}\sinh \left (x_R \sqrt {\eta ^2 +2r} \ \right ) + e^ {-(b+x_R) \eta}\sinh \left ((b - x_R) \sqrt {\eta ^2 +2r}  \ \right ) \right ]}.
\end{equation}
Note that, unlike the case  $\eta =0,$  one has $T(x,r) \neq T(b-x, r),$ namely $T(x,r)$ is not symmetric with respect to $x=b/2.$ \par
The results concerning the minimization of the FET are in part similar to those obtained
for $\eta =0; $ therefore, for the sake of  brevity we limit ourselves to report only few examples in tabular form, without graphical outputs. \par
In the following tables we report the values of  $r_m(x), \ m(x)$ and $T(x,0)= E [\tau (x,0)]$  as functions of $x,$ for $b=1$ and several values of $\eta$ and $x_R \in (0,1).$
\par\noindent
As for $T(x,0),$ one has (see \cite{abundo:MCAP2025}, Eq. (2.44)):
\begin{equation} \label{TTx0}
T(x,0) = \frac 1 \eta \left [\frac {b(1- \exp (-2 \eta x )} {1 - \exp (-2 \eta b )} -x \right ].
\end{equation}  \par\noindent
In Table \ref{tabdriftFETmu1xR02} we report the values of  $r_m(x), \ m(x)$ and $T(x,0)= E [\tau (x,0)]$  as functions of $x,$ from $0$ to $1$ with step $0.1 ,$ for $b=1, \ \eta =1,$ and $x_R =0.2.$
Of course, the numbers in the third column are not greater than those in the fourth one.
\begin{table}[!h]
\begin{center}
\begin{tabular}{cccc}
\hline
 $x$ & $r_{m}(x)$ & $m(x)$ & $T(x,0)$ \\
\hline
0      & 0      & 0       & 0 \\
0.1    & 14.242 & 0.107 & 0.109 \\
0.2    & 30.681 & 0.156 & 0.181 \\
0.3    & 39.731 & 0.122  & 0.221\\
0.4    & 43.269 & 0.1775  & 0.231\\
0.45    & 43.621 & 0.1780 & 0.236\\
0.5    & 43.186 & 0.1778 & 0.231\\
0.6    & 39.642 & 0.1746 & 0.208\\
0.7    & 30.222 & 0.1646 & 0.1713\\
0.8    & 0 & 0.123 & 0.123\\
0.9    & 0 & 0.065 & 0.065 \\
1    & 0 & 0 & 0\\
\hline
\end{tabular}
\end{center}
\caption{For drifted BM with resetting, the table reports the values of $r_{m}(x), \ m(x)=T(x,r_{m}(x))$ and $T(x,0)$ numerically obtained for some  $x \in [0,1],$ with $b=1, \ \eta =1,$ and $x_R =0.2 \ .$ Here, $T(x,r)$ is the expected FET.}
\label{tabdriftFETmu1xR02}
\end{table}
We see that $r_m(x)$ becomes zero for $x > 0.7 \ .$ Moreover, $r_m(x)$ and $m(x)$ are both increasing for $x \in (0, 0.45),$ and decreasing for $x \in (0.45, 1),$
so  at $x=0.45$ there is the maximum; a similar situation was detected in the preceding case  $\eta =0 ,$ but the point of maximum was at $x = 0.5 \ .$ \par
We carried on calculations of the above quantities for $b=1, \ \eta =1,$ and several values of $x_R .$ Actually, for $x_R > 0.2$ we  obtained $r_{m}(x)  =0 ,$ for  all $x.$ \par
Another example is shown in Table \ref{tabdriftFETmu-1xR02}, concerning  the case when $b=1, \ \eta =-1$ and $x_R=0.2 \ ; $
$T(x,0)$ is given again by \eqref{TTx0}.
\begin{table}[!h]
\begin{center}
\begin{tabular}{cccc}
\hline
 $x$ & $r_{m}(x)$ & $m(x)$ & $T(x,0)$ \\
\hline
0      & 0      & 0       & 0 \\
0.1    & 11.168 & 0.059 & 0.065 \\
0.2    & 24.277 & 0.095 & 0.123 \\
0.3    & 36.140 & 0.112  & 0.171\\
0.4    & 43.998 & 0.118  & 0.208\\
0.5    & 45.857 & 0.120 & 0.231\\
0.6    & 45.886 & 0.119 & 0.236\\
0.7    & 42.871 & 0.116 & 0.221\\
0.8    & 35.547 & 0.106 & 0.181\\
0.9    & 24.467 & 0.075 & 0.109 \\
1    & 0 & 0 & 0\\
\hline
\end{tabular}
\end{center}
\caption{For drifted BM with resetting, values of $r_{m}(x), \ m(x)=T(x,r_{m}(x))$ and $T(x,0)$ numerically obtained for some  $x \in [0,1],$ with $b=1, \ \eta =-1,$ and $x_R =0.2 \ .$ here, $T(x,r)$ is the expected FET.}
\label{tabdriftFETmu-1xR02}
\end{table}
Note that now both $r_m(x)$ and $m(x)$ attain their maximum at $x=0.6 \ .$
\par
As before, taking fixed  $b=1$ and $ \eta =-1,$ we  carried on calculations for several values of $x_R \in (0,1) .$
 Actually, we  detected a value $\bar x_R \in (0.23, 0.24)$
such that, when $x_R > \bar x_R ,$ \ $r_{m}(x)$ becomes zero for  all $x;$ as in the case of $\eta =0,$
this can be explained, by considering that, as $x_R$ becomes close enough to $1/2 ,$ also small non-zero values of the resetting rate $r$ imply large values of
the expected FET. In fact, if the process is reset  to a position close to $1/2,$ it takes more time to reach one of the ends of the interval $(0, 1),$ hence the minimum of the expected FET is obtained at $r=0.$ \par
Of course, this kind of behavior was observed for any values of $b$ and $\eta ,$  when $x_R$ is close to $b/2.$ \par\noindent

\subsection{The case of OU with resetting: optimization of the FPT }
Let $X(t)$ be OU process driven by the SDE \eqref{OUSDE} and ${\cal X}(t)$ the corresponding process with resetting.
We limit ourselves to study minimization of the expected FPT, since formula \eqref{meanFETOU} for the mean FET of OU with resetting is rather computationally complicated.
\bigskip

\noindent {\bf (i) The case when $x=x_R$ } \par\noindent
Now, formula \eqref{EFPTOUx=XRreset} holds for the expected FPT, namely:
\begin{equation}
T(x,r):=  E[\tau (x,r)] = \frac 1 r \ \left [ \frac {D_ {- r / \mu} ( 0)}  {D_ {- r / \mu} \left (x \sqrt {2 \mu / \sigma ^2} \right )}  \ e^  {\frac {- \mu x^2 } {2 \sigma ^2} } -1 \right ].
  \end{equation}
As in the case of BM with resetting,  for fixed starting point $x = x_R >0,$ the quantity $T(x,r),$   as a function of $r,$ attains its unique global minimum at a value
$
r_m(x)= arg \left ( \min _ {r \ge 0 } T (x,r) \right ).
$
As before, our goal is to find $r_{m}(x)$ and the
minimum expected FPT, $m(x)= T ( x,r_{m}(x)) ,$ for fixed $x = x_R  >0.$
Actually, in the computations we make use of the WolframAlpha software, because the special function $D_ \nu$ is involved. \par
\begin{figure}[!h]
\centering
\includegraphics[height=0.35 \textheight]{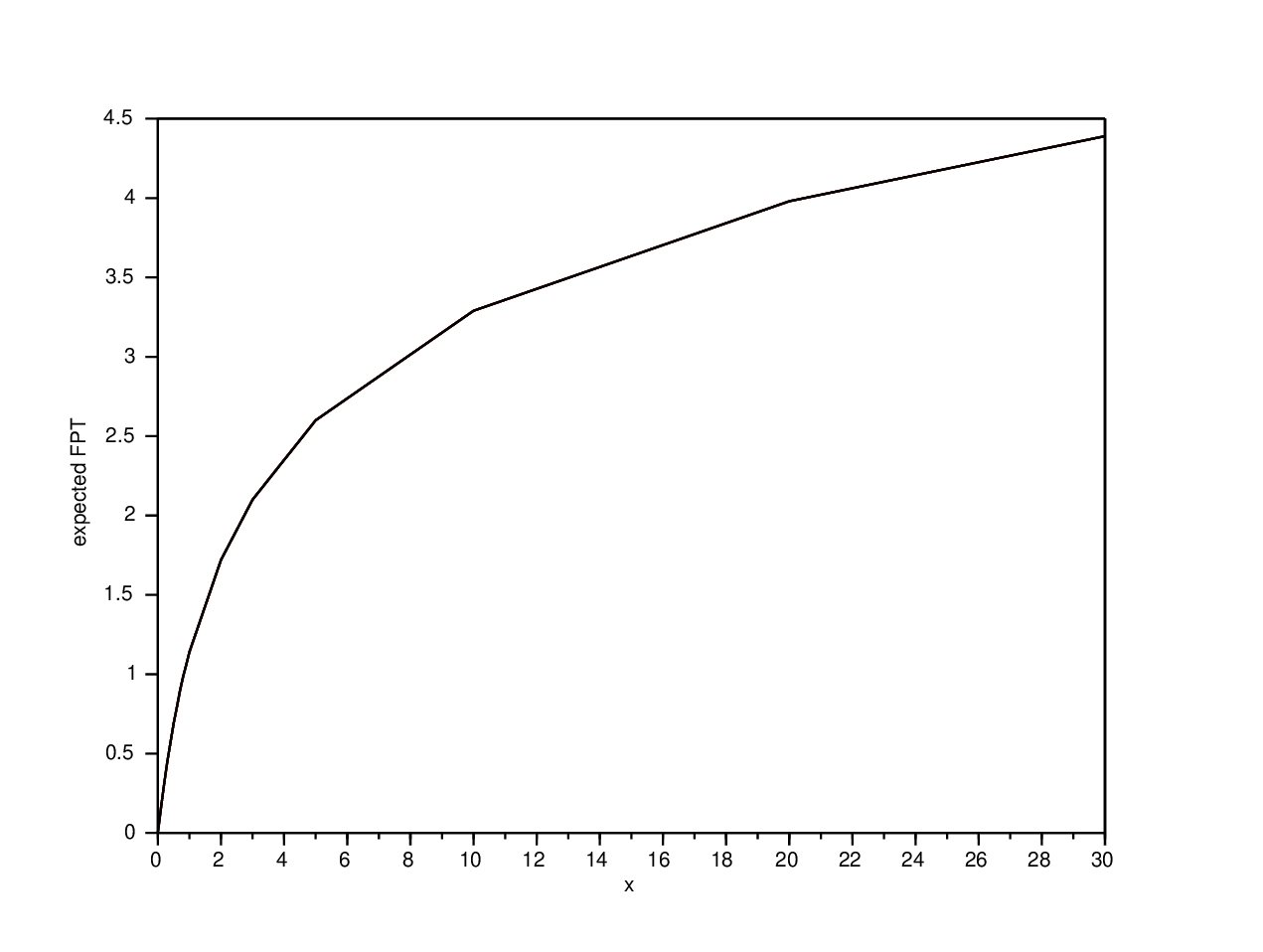}
\caption{Graph of the expected FPT, $T(x,0),$ of OU process without resetting with $\mu = \sigma =1,$  as function of $x=x_R.$
}
\label{FPTOUr=0}
\end{figure}

In the no-resetting case $(r=0)$ with $\mu= \sigma =1,$ we report  in Figure \ref{FPTOUr=0}, the graph of the expected FPT, $T(x,0),$ as a function of $x= x_R.$
We see that $T(x,0)$ behaves as a concave, increasing function of $x.$


Now, we go to consider the case $r \neq 0.$
In Table \ref{tabOUr}, for $\mu = \sigma =1,$ and some values of $x=x_R >0 $ (in the first column), we report in the second one the value $r_{m}(x)$ which minimizes $T(x,r), $ in the third column the optimal value $m(x) = \min _ {r \ge 0 } T (x,r)= T(x, r_{m}(x)),$
and in the fourth one the value $T(x,0),$ i.e. the expected FPT in the no-resetting case; all values was  obtained by using the WolframAlpha software.
\par
\begin{figure}
\centering
\includegraphics[height=0.35 \textheight]{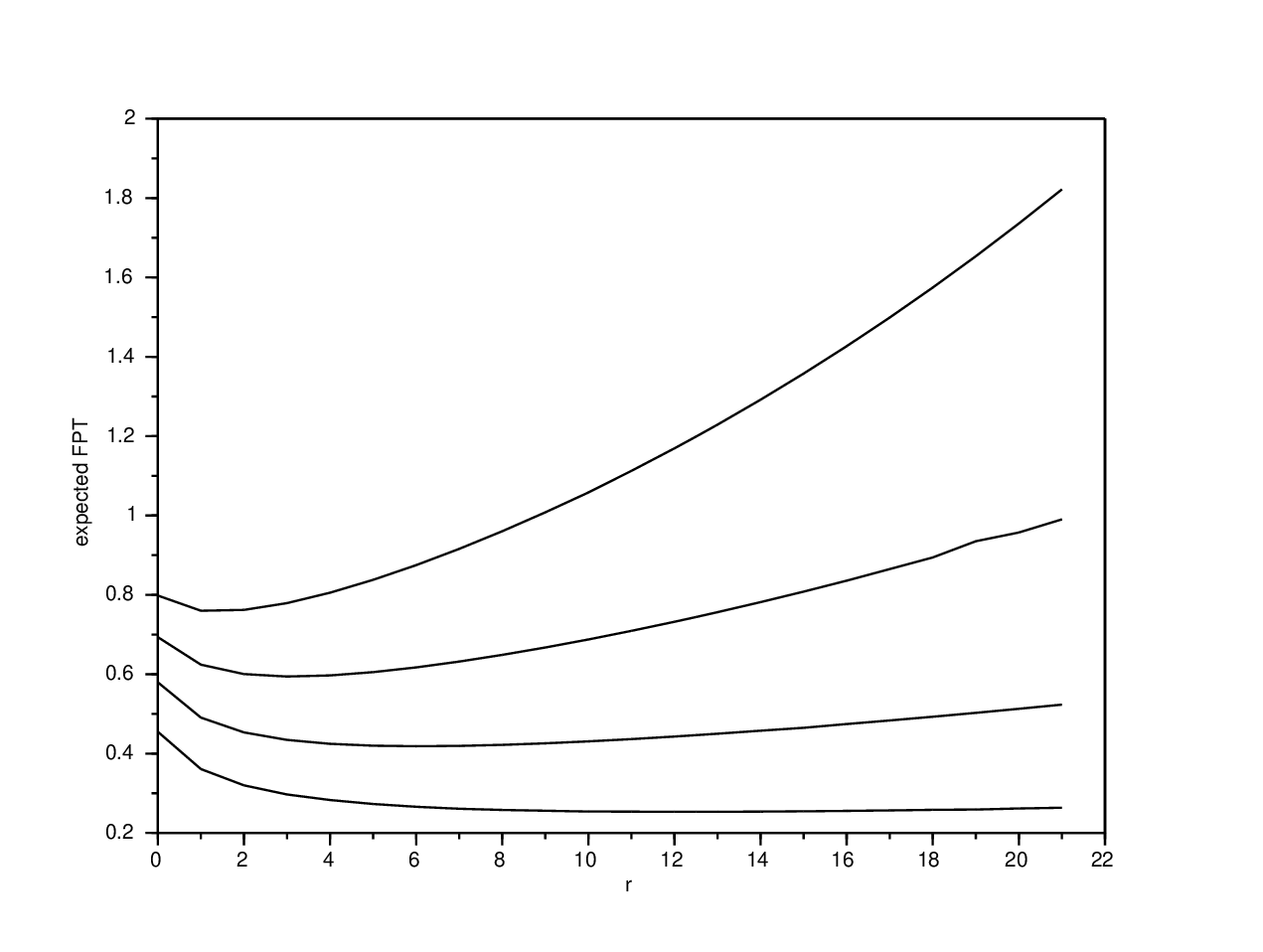}
 \caption{For OU process with resetting and $\mu = \sigma =1,$ the figure shows the graphs of the expected FPT, $T(x,r),$ as a function of $r,$ for the values of $x$ from $0.3$ to $0.6$ with step $0.1 \ ;$
  on the horizontal axes $r,$ on the vertical one $T(x,r).$
 The curves are ordered in this way: the greater the value of $x,$ the higher the corresponding curve  }
 \label{grafOUfour}
\end{figure}

\noindent
In Figure \ref{grafOUfour}, for OU process with resetting and $\mu = \sigma =1,$  we report four examples of the graphs of $T(x,r),$ as a function of $r,$
for the values of $x = x_R$ from $0.3$ to $0.6$ with step $0.1 \ .$
The curves are ordered in this way: the greater the value of $x,$ the higher the corresponding curve. The minimum values $r_m(x)$ of $T(x,r)$ are reported in Table \ref{tabOUr}; as we can see, $r_m(x)$ shifts towards the left, as $x$ increases (unlike the case of BM with resetting), while $m(x)$ shifts upward.
Of course,  we set $T(0,r) =0,$ for any $r.$ As we see,
$r_{m}(x)$ decreases, as a function of $x,$ while $m(x)$ and $T(x,0)$
increase.


\par\noindent
Note that,
for $x < 0.8$  all the numbers in the third column of Table \ref{tabOUr} are less  than the corresponding ones  in the fourth column, showing that  optimal resetting really reduces the mean FPT (with respect to the no-resetting case); instead
for $x \ge 0.8$  the values in the second column become zero, hence those in third column become equal to the corresponding ones  in the fourth column. This can be explained by considering that, if $x = x_R$ is
large enough, resetting to the position $x_R$ at a positive rate, instead of speeding up, slow down  the first passage through zero. \par
\begin{figure}[!h]
\centering
\includegraphics[height=0.35 \textheight]{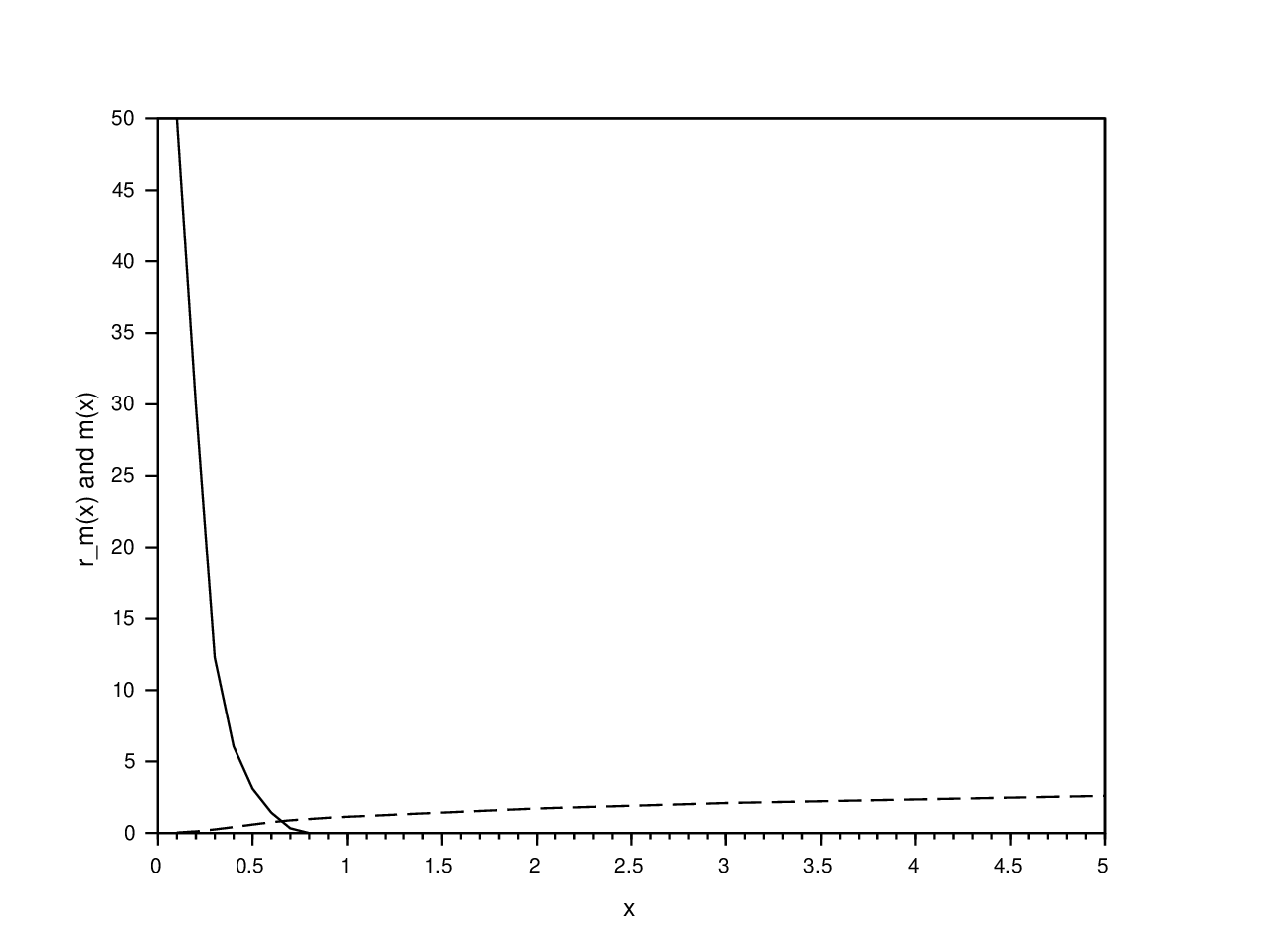}

\caption{For OU process with resetting, and $\mu = \sigma =1 ,$ the figure shows the  graphs of $r_m(x)$ (solid line) and $m(x)$ (dashed line), for $x \in [0.1, 5].$
}
\label{EFPTOU11_x}
\end{figure}

\noindent
In Figure \ref{EFPTOU11_x} we report the graphs of $r_m(x)$ (solid line) and $m(x)$ (dashed line), as functions of $x,$ for OU process with resetting and  $\mu = \sigma =1$ (the values are extracted  from Table \ref{tabOUr}).  \par
\begin{figure}[!h]
\centering
\includegraphics[height=0.35 \textheight]{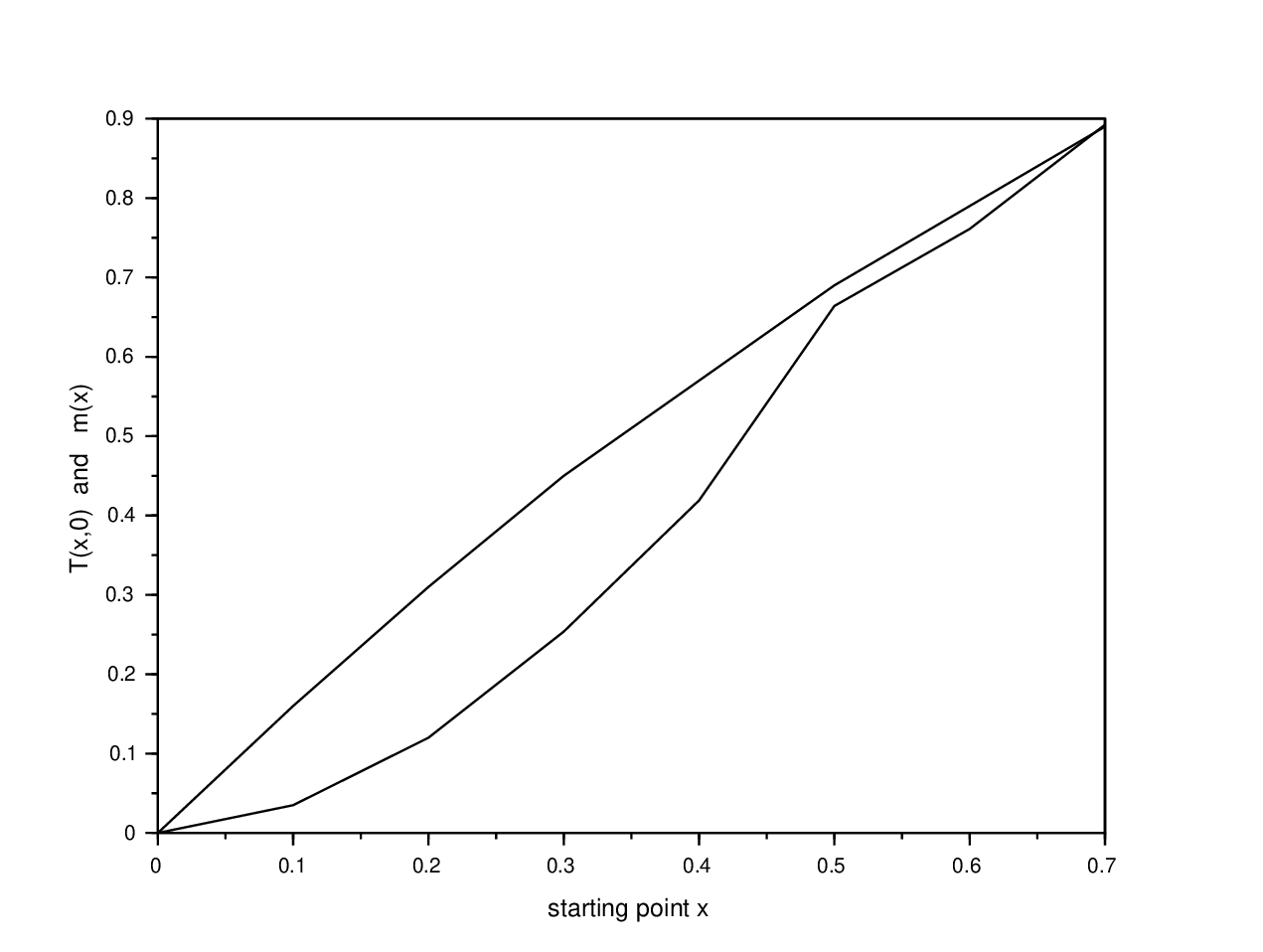}

\caption{For OU process with resetting, and $\mu = \sigma =1 ,$ the figure shows the comparison between the graph of $T(x,0),$ i.e. the expected FPT in the no-resetting case  (higher curve) and $m(x)$ (lower curve), for $x \in [0.1, 0.7].$
}
\label{FPTOUcomparison}
\end{figure}

\noindent
In Figure \ref{FPTOUcomparison} we show the comparison between the graph of $T(x,0)$ (higher curve) and that of $m(x)$ (lower curve), for $x \in [0.1, 0.7];$ we have not reported the graphs for $x \ge 0.8,$ since for such values
the curves of $T(x,0)$ and
$m(x)$ coincide, that is, $T(x,0)\equiv m(x)$ (see Table \ref{tabOUr}).

\begin{table}[!h]
\begin{center}
\begin{tabular}{cccc}
\hline
 $x$ & $r_{m}(x)$ & $m(x)$ &  $T(x,0)$ \\
\hline
0.     & 0. & 0. & 0. \\
0.1    & 49.99 & 0.03 & 0.16 \\
0.2    & 29.97 & 0.11 & 0.31 \\
0.3    & 12.29 & 0.25 & 0.54 \\
0.4    & 6.05 & 0.41 & 0.57 \\
0.5    & 3.10 & 0.59 & 0.69 \\
0.6    & 1.42 & 0.75 & 0.79 \\
0.7    & 0.33 & 0.89 & 0.891 \\
0.8    & 0. & 0.98 & 0.98 \\
0.9    & 0. & 1.07 & 1.07 \\
1.     & 0. & 1.14 & 1.14  \\
2.     & 0. & 1.72 & 1.72 \\
3.     & 0. & 2.10 & 2.10  \\
5.     & 0. & 2.60 & 2.60 \\
10.     & 0. & 3.29  & 3.29 \\
20.     & 0. & 3.98  & 3.98\\
30.     & 0. & 4.38 & 4.38  \\
\hline
\end{tabular}
\end{center}
\caption{For OU process with resetting, and $\sigma = \mu = 1,$ the table reports the values of $r_{m}(x), \ m(x)=T(x,r_{m}(x)),$  and $T(x,0),$ for some values of $x =x_R >0.$ Here, $T(x,r)$ is the expected FPT. }
 \label{tabOUr}
\end{table}

\begin{table}[!h]
\begin{center}
\begin{tabular}{cccc}
\hline
 $x$ & $r_{m}(x)$ & $m(x)$ &  $T(x,0)$ \\
\hline
0.     & 0. & 0. & 0. \\
0.1    & 29.99 & 0.03 & 0.54 \\
0.2    & 29.98 & 0.12 & 1.07 \\
0.3    & 13.93 & 0.27 & 1.58 \\
0.4    & 7.76 & 0.48 & 2.08 \\
0.5    & 4.90 & 0.75 & 2.56 \\
0.6    & 3.35 & 1.07 & 3.03 \\
0.7    & 2.41 & 1.43 & 3.47 \\
0.75   & 2.07 & 1.63 & 3.69 \\
0.8    & 1.80 & 1.85 & 3.91 \\
0.9    & 1.38 & 2.30 & 4.33 \\
1.     & 1.08 &  2.78 & 4.74 \\
1.5    & 0.36 &  5.49 & 6.64 \\
2.     & 0.10 &   8.06 & 8.30 \\
3.     & 0. & 11.09 & 11.09 \\
5.     & 0. & 15.26 & 15.26 \\
10.    & 0. & 21.73 & 21.73 \\
20.    & 0. & 28.66 & 28.66 \\
\hline
\end{tabular}
\end{center}
\caption{For OU process with resetting, and $\sigma =1, \ \mu =0.1 ,$ the table reports the values of $r_{m}(x), \ m(x)=T(x,r_{m}(x)),$  and $T(x,0),$ for some values of $x =x_R >0.$ Here, $T(x,r)$ is the expected FPT. }
 \label{tabOUrmu=01}
\end{table}
In Table \ref{tabOUrmu=01}, instead, we refer to OU process with resetting, and $\sigma =1, \ \mu =0.1 ;$ we report the values of $r_{m}(x), \ m(x)=T(x,r_{m}(x)),$  and $T(x,0),$ for some values of $x =x_R >0.$ \par\noindent
As it must be, all the numbers in  the third column are not greater than those in  the forth column; in fact, the third column reports the minimum value of the expected FPT, attained at the optimal value $r = r_m(x),$ while
the forth column reports the  expected FPT in the case without resetting $(r=0).$ \par\noindent
By comparing Table \ref{tabOUr} (which refers to the case $\mu =1)$ with Table \ref{tabOUrmu=01} (concerning the case $\mu =0.1),$ we  note that, in Table \ref{tabOUr} the minimum values, $m(x),$ of the expected FPT
are smaller than those in Table \ref{tabOUrmu=01}. This is right, because in the first case the drift parameter  $\mu =1 $ is greater than that of the second case $(\mu =0.1),$ so the process is able to reach more quickly the origin.
\par\noindent
In Table \ref{tabOUrmu=005}, for OU process with resetting and $\sigma =1, \ \mu = 0.05,$ we report again the values of $r_m(x), \  m(x)=T(x,r_{m}(x)),$  and $T(x,0),$
for some values of $x =x_R >0.$ Now, the minimum values, $m(x),$ of the expected FPT
are greater than those in Table \ref{tabOUrmu=01}, regarding the case $\mu = 0.1 \ .$ Moreover,  we can observe that for $x$ near enough to $1,$ the values are rather similar to those in Table \ref{tab1}, concerning BM with resetting; in fact, $\mu = 0.05$ is small enough, so the situation is close to the case $\mu =0,$ that is BM with resetting. \par\noindent
Unfortunately, we were unable to study the cases for values of $0 < \mu < 0.05 ,$ because for such values the software WolframAlpha failed to make minimization of the function $T(x,r).$
\par
In all cases, for  OU process with resetting ${\cal X}(t),$ with $x=x_R ,$ the minimum expected FPT through the origin appears to be smaller than  that regarding BM with resetting,
due to the presence of the drift coefficient $\mu (x)= - \mu x;$ moreover, as for BM with resetting in the case  $x=x_R, \ m(x)$ increases, while $r_m(x)$ decreases (cf. Table \ref{tab1}).


\begin{table}[!h]
\begin{center}
\begin{tabular}{cccc}
\hline
 $x$ & $r_{m}(x)$ & $m(x)$ &  $T(x,0)$ \\
\hline
0.     & 0. & 0. & 0. \\
0.1    & 9.9999 & 0.05626& 0.7726 \\
0.2    & 9.9997 & 0.14407& 1.5271 \\
0.3    & 9.9996 & 0.28103& 2.2642 \\
0.4    & 7.8495 & 0.4901& 2.9844 \\
0.5    & 4.9992 & 0.7622& 3.6885 \\
0.6    & 3.4399 & 1.0915& 4.3770 \\
0.7    & 2.5038 & 1.4757& 5.0504 \\
0.8    & 1.8961 & 1.9127& 5.7092 \\
0.9    & 1.4793 & 2.3997& 6.3539 \\
1    & 1.18107 & 2.9338& 6.9851 \\
1.5    & 0.47262 & 6.187& 9.9526 \\
2    & 0.2212 & 10.0305& 12.6418\\
3    & 0.0306 & 17.1461& 17.3407\\
5    & 0. & 24.7728& 24.7728\\
10    & 0.& 37.0560& 37.0560\\
20    & 0.& 50.8449& 50.8449\\
\hline
\end{tabular}
\end{center}
\caption{For OU process with resetting, and $\sigma =1, \mu = 0.05, \ $ the table reports the values of $r_{m}(x), \ m(x)=T(x,r_{m}(x)),$  and $T(x,0),$ for some values of $x =x_R >0.$ Here, $T(x,r)$ is the expected FPT. }
 \label{tabOUrmu=005}
\end{table}
\bigskip
\newpage
\noindent {\bf (ii) The case $x \neq x_R$} \par\noindent
Now, formula \eqref{EFPTOUreset} holds for the expected FPT through the origin, namely:
 \begin{equation}
 T(x,r) = E[\tau (x,r)] = \frac 1 r \left [ \left (1- e^  {\frac {\mu x^2 } {2 \sigma ^2} } \  \frac {D_ {- r / \mu} \left (x \sqrt {2 \mu / \sigma ^2} \right )} {D_ {- r / \mu} ( 0)}  \right ) \cdot
 \frac 1 {e^  {\frac {\mu x^2 } {2 \sigma ^2} } \  \frac {D_ {- r / \mu} \left (x_R \sqrt {2 \mu / \sigma ^2} \right )} {D_ {- r / \mu} ( 0)}    }    \right ].
 \end{equation}
As in the previous case,  for fixed starting point $x  >0, \ T(x,r),$   as a function of $r,$ attains its unique global minimum at the value
$
r_m(x)= arg \left ( \min _ {r \ge 0 } T (x,r) \right ).
$
We will search for the
minimum expected FPT, $m(x)= T ( x,r_{m}(x)) ,$ for fixed $x  >0.$
Also now, we will make use of the WolframAlpha software. \par
In Table \ref{tabOUxR1}, for $\mu = \sigma =1, \ x_R=1 $ and some values of $x$ (in the first column), we report in the second one the value $r_{m}(x)$ which minimizes $T(x,r), $ in the third column the optimal
value $m(x) = \min _ {r \ge 0 } T (x,r)= T(x, r_{m}(x)),$
and in the fourth one the value $T(x,0),$ i.e. the expected FPT in the no-resetting case (all the values are approximated up to the third decimal digit).
\begin{table}[!h]
\begin{center}
\begin{tabular}{cccc}
\hline
 $x$ & $r_{m}(x)$ & $m(x)$ &  $T(x,0)$ \\
\hline
0.     & 0. & 0. & 0. \\
0.1    & 0. & 0.167& 0.167 \\
0.2    & 0. & 0.318& 0.318 \\
0.3    & 0 & 0.455& 0.455 \\
0.4    & 0. & 0.579& 0.579 \\
0.5    & 0 & 0.693& 0.693 \\
1    & 0. & 1.147& 1.147 \\
1.5    & 0. & 1.475& 1.475 \\
2    & 0.388 & 1.714& 1.728\\
2.5 & 0.888 & 1.844& 1.933 \\
3    & 1.195 & 1.916& 2.105\\
5    & 1.692 & 2.016& 2.599\\
7    & 1.843 & 2.041& 2.931\\
10    & 1.928& 2.054& 3.284\\
15    & 1.977& 2.060& 3.687\\
20    & 1.996& 2.062& 3.974\\
50    & 2.023& 2.064& 4.886\\
\hline
\end{tabular}
\end{center}
\caption{For OU process with resetting, $\sigma = \mu = 1 $ and $x_R=1$  the table reports the values of $r_{m}(x), \ m(x)=T(x,r_{m}(x)),$  and $T(x,0),$ for some values of $x >0.$ Here, $T(x,r)$ is the expected FPT. }
 \label{tabOUxR1}
\end{table}

We see that $m(x) \le T(x,0),$ as it must be; moreover, $r_m(x)$ and $m(x)$ are both increasing functions of $x,$ as in the case of BM with resetting with $x \neq x_R$ (see Table \ref{tab2}). \par\noindent
Unlike the case of BM, we observe that only for $x \ge 2$ it results $r_m(x) \neq 0;$ this can be explained, because if $x$ is close enough to $x_R=1,$ resetting implies wasting time in reaching zero. Actually, the situation was different for BM; now there is the effect of
the drift coefficient $\mu (x)= - \mu x.$
Furthermore, for large $x>0, \ r_m(x)$ and $m(x)$ approach two asymptotic values (precisely, for $x=50$ they take on the values $2.023 $ and $2.064,$ respectively). \par
\begin{figure}[!h]
\centering
\includegraphics[height=0.35 \textheight]{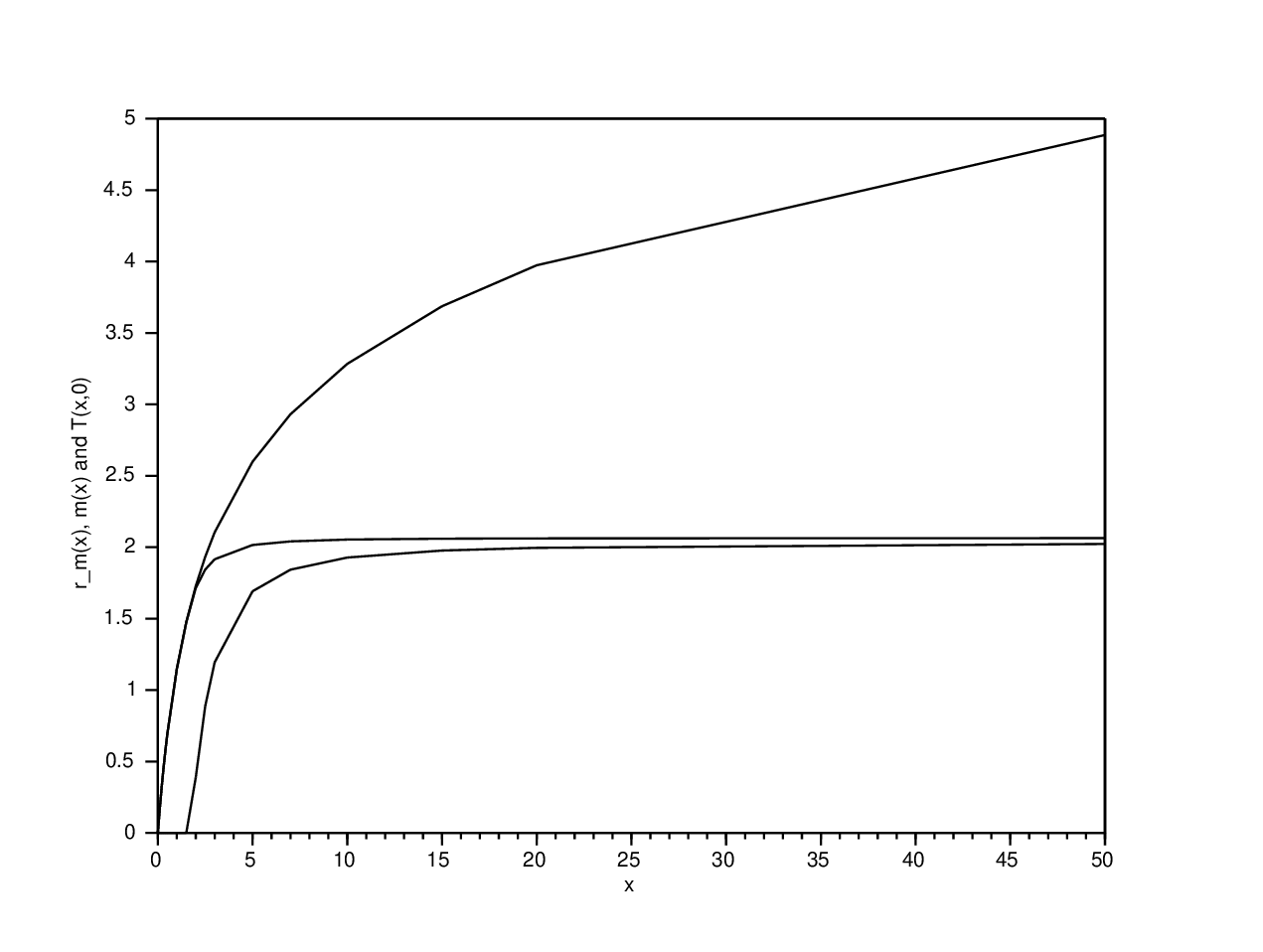}

\caption{For OU process with $\mu = \sigma =1$ and resetting position $x_R=1,$  the figure shows the comparison between the graphs of $r(x)$ (lower curve), $m(x)$ (middle curve) and $T(x,0)$ (higher curve), for $x \in [0,50].$ Note that for $x<2, \ r(x)$ is zero, so $m(x)$ coincides with $T(x,0).$
}
\label{compOUxR1}
\end{figure}

\noindent
Taking even larger values of $x, \ T(x,0)$ seems to be indefinitely increasing with $x$ (see Figure \ref{compOUxR1}).

\par
In Table \ref{tabOUxR2}, for $\mu = \sigma =1, \ x_R=2 $ and some values of $x$ (in the first column), we report in the second one the value $r_{m}(x)$ which minimizes $T(x,r), $ in the third column the optimal
value $m(x) = \min _ {r \ge 0 } T (x,r)= T(x, r_{m}(x)),$
and in the fourth one the value $T(x,0)$  (all the values, except that in third column  corresponding to $x=8,$  are approximated up to the third decimal digit).
\begin{table}[!h]
\begin{center}
\begin{tabular}{cccc}
\hline
 $x$ & $r_{m}(x)$ & $m(x)$ &  $T(x,0)$ \\
\hline
0.     & 0. & 0. & 0. \\
0.1    & 0. & 0.167& 0.167 \\
0.2    & 0. & 0.318& 0.318 \\
0.3    & 0 & 0.455& 0.455 \\
0.4    & 0. & 0.579& 0.579 \\
0.5    & 0 & 0.693& 0.693 \\
1    & 0. & 1.147& 1.147 \\
1.5    & 0. & 1.475& 1.475 \\
2    & 0. & 1.728& 1.728\\
2.5 & 0. & 1.844& 1.933 \\
3    & 0. & 2.106& 2.106\\
5    & 0. & 2.600& 2.600\\
7    & 0. & 2.932& 2.932\\
8    & 0.008 & 3.06500& 3.06505\\
10    & 0.193& 3.255& 3.286\\
15    & 0.421& 3.488& 3.689\\
20    & 0.527& 3.594& 3.976\\
50    & 0.849& 3.876& 5.580\\
\hline
\end{tabular}
\end{center}
\caption{For OU process with resetting, $\sigma = \mu = 1 $ and $x_R=2$  the table reports the values of $r_{m}(x), \ m(x)=T(x,r_{m}(x)),$  and $T(x,0),$ for some values of $x >0.$ Here, $T(x,r)$ is the expected FPT. }
 \label{tabOUxR2}
\end{table}
We see that $m(x) \le T(x,0),$ and  $r_m(x)$ and $m(x)$ are both increasing functions of $x,$ as in the previous case (see Table \ref{tabOUxR1}). Now,   only for $x \ge 8$ it results $r_m(x) \neq 0;$ this can be explained
by the same considerations done in the case of drifted BM with resetting and $x= x_R$ (see pg. 17)
\par\noindent
For large $x>0, \ r_m(x)$ approaches the asymptotic value $0.849  ,$ while $m(x)$ approaches the  value $3.876 ;$ taking even larger values of $x, \ T(x,0)$ appears to be indefinitely increasing with $x.$

\section {Conclusions and final Remarks}
This study investigated the problem of minimizing the expected first-passage time (FPT) through zero and the expected first-exit time (FET) from an interval $(0,b),$ for a diffusion process with Poissonian resetting, denoted as ${\cal X}(t)$. This process is derived from a one-dimensional, time-homogeneous diffusion $X(t)$. We specifically examined cases where $X(t)$ is either a Brownian motion (BM) with drift or an Ornstein-Uhlenbeck (OU) process.
\par\noindent
Our motivation stemmed from the frequent need in fields like biological diffusion models and mathematical finance to identify the optimal resetting rate $r$ that minimizes the expected FPT (for the one-boundary case) or FET (for the two-boundary case), thereby expediting first-passage. To achieve this optimization, we first established a relation between the Laplace transforms of the FPT (or FET) for ${\cal X}(t)$ and $X(t)$. From this, we derived the first and second moments of the FPT and FET for ${\cal X}(t)$ in terms of the Laplace transforms of the FPT and FET for $X(t)$. This approach proved more computationally tractable than directly solving the differential problems for these moments (see Eqs. \eqref{eqmeanFPT}, \eqref{eqsecondFPT}), which can be rather complex.
\par
For drifted Brownian motion with resetting, we explored the minimization of both the expected FPT and the expected FET. For the OU process with resetting, we focused solely on minimizing the expected FPT due to the increased complexity of the expected FET formula. \bigskip

Regarding the minimization of the expected FPT for undrifted BM with resetting, our theoretical and numerical analyses revealed that for a fixed reset position $x_R > 0$, the optimal resetting rate $r_m(x)$ (which minimizes the expected FPT) is an increasing function of $x > 0$. This $r_m(x)$ falls within an interval $(\alpha(x), \beta(x))$, where both $\alpha(x)$ and $\beta(x)$ are decreasing functions of $x_R$. Furthermore, for a fixed $x_R$, $r_m(x)$ exhibited a jump discontinuity at $x=0$. Our study also demonstrated that for $x > 0$, the minimum expected FPT, $m(x),$ increases from $m(0)=0$ to $m(\infty) = \frac{1}{2} e^2 x_R^2$.
\bigskip

In contrast, the expected FET of Brownian motion with resetting presented a more intricate scenario for the minimum $m(x)$ and its corresponding optimal resetting rate $r_m(x),$ when the reset position $x_R \in (0,b)$ was fixed. Numerical computations allowed us to elucidate their qualitative behaviors. For instance, with $x_R = 0.2$, both $r_m(x)$ and $m(x)$ attained their maximum at $x=1/2$. As $x$ increased, the expected FET curves, as functions of $r>0,$ became progressively higher, and the value of $r_m(x)$ shifted further to the right. However, we observed that the qualitative behaviors of $r_m(x)$ and $m(x)$ remained independent of $x_R,$ unless $x_R$ was  close enough to $1/2$. We identified a transition value $\bar x_R = 0.295$, at which a change in behavior occurred. Specifically, for $x_R = \bar x_R$, the maximum value of $r_m(x)$ became zero for all $x$, and $m(x)$ coincided with $x(1-x)$, that is the expected FET in the absence of resetting ($r=0$). Although our presented results used $b=1$, this phenomenon holds for any value of $b$, with the behavioral change observed when $\bar x _R$ is sufficiently close to $b/2$.\bigskip

For the minimization of the expected FPT and FET of drifted BM with resetting, we found results consistent with those observed for the zero-drift case.
\bigskip

Regarding the minimization of the expected FPT of the OU process with resetting, we set $\sigma=1$ and compared the results for $\mu = 0.05$, $\mu = 0.1$, and $\mu = 1$. We noted that the minimum values, $m(x)$, of the expected FPT  were greater in the first two cases (smaller drift parameters) than in the third case, allowing the process to reach the origin more slowly. \par\noindent
In the case $x=x_R,$ the minimum expected FPT, $m(x),$ appeared smaller than that for (undrifted) BM with resetting, especially for large $x.$ Conversely, when the BM with resetting exhibited a sufficiently large negative drift, the opposite effect was observes (cf. e.g. the results for BM with drift $\eta =-1$
and those for OU with  $\mu =1.$
\bigskip

In conclusion, this work confirms that a resetting mechanism substantially expedites the first-passage of a diffusion process through one or two boundaries. Therefore, investigating the optimal resetting rate $r$ that minimizes the mean FPT and FET is a worthwhile endeavor. While this study focused on drifted BM and OU processes, future research could explore other diffusion processes with resetting. However, in such cases, the absence of closed-form expressions for the expected FPT and FET would necessitate the use of Monte Carlo simulations. Finally, it would be valuable to investigate the simultaneous optimization of the
mean and the variance of the FPT (or FET) to identify a resetting rate $r$ that minimizes both quantities.


\section*{Acknowledgments}

The author belongs to GNAMPA, the Italian National Research Group of INdAM; he also acknowledges
the MUR Excellence Department Project MatMod@TOV awarded to the Department of Mathematics, University of Rome Tor Vergata, CUP E83C23000330006 .



\end{document}